\documentclass[12pt]{amsart}
\title[Affine Markov trace]{Markov trace on a tower of affine Temperley-Lieb algebras of type $\tilde{A}$}
\author{Sadek AL HARBAT}
\address{LAMFA, Universit\'e de Picardie - Jules Verne} 
\email{sadikharbat@math.univ-paris-diderot.fr}

\usepackage{etex} 
\usepackage[latin1]{inputenc}	
\usepackage[T1]{fontenc}
\usepackage[english]{babel}
\usepackage{amsmath,amssymb}
\usepackage{wasysym}
\usepackage{stmaryrd}
\usepackage[table]{xcolor}
\usepackage{color}
\usepackage{dsfont}\let\mathbb\mathds
\usepackage[all]{xy}
\usepackage{graphicx}
\usepackage{mathenv}
\usepackage{lastpage}
\usepackage{fancyhdr}	
\usepackage{subfigure}
\usepackage{geometry}
\usepackage[version=3]{mhchem}
\usepackage[force,almostfull]{textcomp}
\usepackage{array}
\usepackage{lipsum}		
\usepackage{eso-pic}
\usepackage{multirow}
\usepackage{fancybox}
\usepackage{cite}
\usepackage{amsthm}
\usepackage{listings}
\usepackage{lscape}
\usepackage{multicol}		
\usepackage{setspace}
\usepackage{times}
\usepackage{eurosym}
\usepackage{latexsym}
\usepackage{ulem}
\usepackage{lmodern}
\usepackage{colortbl}
\usepackage{rotating}
\usepackage{type1cm}
\usepackage{lettrine}
\usepackage{ragged2e}
\usepackage{mathrsfs}
\usepackage{scrtime}
\usepackage{enumitem}

\usepackage{longtable}
\usepackage{url}
\usepackage{nomentbl}
\usepackage{nomencl}
\usepackage{hyperref}
\hypersetup{
pdfpagemode=UseOutlines,      
pdfstartview=Fit,             
pdffitwindow=true,            
pdfpagelayout=TwoColumnsRight,
pdftoolbar=true,              
pdfmenubar=true,              
bookmarksopen=false,          
bookmarksnumbered=true,       
colorlinks=true,              
pdfauthor={Ton nom},          
pdftitle={Titre PDF},         
pdfcreator=PDFLaTeX,          %
pdfproducer=PDFLaTeX,         %
linkcolor=blue,               
urlcolor=blue,                
anchorcolor=black,            
citecolor=blue,               
frenchlinks=true,             
pdfborder={0 0 0}             
}
\usepackage{pgf, tikz}

\usetikzlibrary{matrix,positioning}
\usetikzlibrary{arrows,decorations.markings}

\makeatletter
\newlength\@tempdim@x
\newlength\@tempdim@y
\newcommand\AtUpperLeftCorner[3]{%
\begingroup
\@tempdim@x=0cm
\@tempdim@y=\paperheight
\advance\@tempdim@x#1
\advance\@tempdim@y-#2
\put(\LenToUnit{\@tempdim@x},\LenToUnit{\@tempdim@y}){#3}%
\endgroup}
\newcommand\AtUpperRightCorner[3]{%
\begingroup
\@tempdim@x=\paperwidth
\@tempdim@y=\paperheight
\advance\@tempdim@x-#1
\advance\@tempdim@y-#2
\put(\LenToUnit{\@tempdim@x},\LenToUnit{\@tempdim@y}){#3}%
\endgroup}
\newcommand\AtLowerLeftCorner[3]{%
\begingroup
\@tempdim@x=0cm
\@tempdim@y=0cm
\advance\@tempdim@x#1
\advance\@tempdim@y#2
\put(\LenToUnit{\@tempdim@x},\LenToUnit{\@tempdim@y}){#3}%
\endgroup}
\newcommand\AtLowerRightCorner[3]{%
\begingroup
\@tempdim@x=\paperwidth
\@tempdim@y=0cm
\advance\@tempdim@x-#1
\advance\@tempdim@y#2
\put(\LenToUnit{\@tempdim@x},\LenToUnit{\@tempdim@y}){#3}%
\endgroup}
\makeatother
\newtheorem{theoreme}{Theorem}[section]

\newtheorem{definition}[theoreme]{Definition}
\newtheorem{proposition}[theoreme]{Proposition}
\newtheorem{lemme}[theoreme]{Lemma}
\newtheorem{corollaire}[theoreme]{Corollary}

\newtheorem{remarque}[theoreme]{Remark}

\newenvironment{demo}{\begin{proof}}{\end{proof}}

\addto\captionsfrenchb{}
\addto\captionsfrenchb{}

    \newlength{\myarrowsize} 
    \newlength{\myoldlinewidth}

    \pgfarrowsdeclare{myto}{myto}{
        \pgfsetdash{}{0pt} 
        \pgfsetbeveljoin 
        \pgfsetroundcap 
        \setlength{\myarrowsize}{0.5pt}
        \addtolength{\myarrowsize}{.5\pgflinewidth}
        \pgfarrowsleftextend{-4\myarrowsize-.5\pgflinewidth} 
        \pgfarrowsrightextend{.7\pgflinewidth}
    }{
        \setlength{\myarrowsize}{0.4pt} 
        \addtolength{\myarrowsize}{.3\pgflinewidth}  
        \setlength{\myoldlinewidth}{\pgflinewidth}
        \pgfsetroundjoin
        \pgfsetlinewidth{0.0001pt}
        \pgfpathmoveto{\pgfpoint{0.43\myarrowsize}{0}}
        \pgfpatharc{0}{70}{0.14\myarrowsize}
        \pgfpatharc{-80}{-169.5}{4\myarrowsize}
        \pgfpatharc{150}{189}{0.95\myarrowsize and 0.95\myarrowsize}
        \pgfpatharc{0}{-40}{15\myarrowsize}
        \pgfpathmoveto{\pgfpoint{0.43\myarrowsize}{0}}
        \pgfpatharc{0}{-70}{0.14\myarrowsize}
        \pgfpatharc{80}{169.5}{4\myarrowsize}
        \pgfpatharc{-150}{-189}{0.95\myarrowsize and 0.95\myarrowsize}
        \pgfpatharc{0}{40}{15\myarrowsize}
        \pgfpathclose
        \pgfsetstrokeopacity{0.25}
        \pgfusepathqfillstroke
    }

    \pgfarrowsdeclare{myonto}{myonto}{
        \pgfsetdash{}{0pt} 
        \pgfsetbeveljoin 
        \pgfsetroundcap 
        \setlength{\myarrowsize}{0.5pt}
        \addtolength{\myarrowsize}{.5\pgflinewidth}
        \pgfarrowsleftextend{-4\myarrowsize-.5\pgflinewidth} 
        \pgfarrowsrightextend{.7\pgflinewidth}
    }{
        \setlength{\myarrowsize}{0.4pt} 
        \addtolength{\myarrowsize}{.3\pgflinewidth}  
        \setlength{\myoldlinewidth}{\pgflinewidth}
        \pgfsetroundjoin
        \pgfsetlinewidth{0.0001pt}
        \pgfpathmoveto{\pgfpoint{0.43\myarrowsize}{0}}
        \pgfpatharc{0}{70}{0.14\myarrowsize}
        \pgfpatharc{-80}{-169.5}{4\myarrowsize}
        \pgfpatharc{150}{189}{0.95\myarrowsize and 0.95\myarrowsize}
        \pgfpatharc{0}{-40}{15\myarrowsize}
        \pgfpathmoveto{\pgfpoint{-7\myarrowsize}{0}}
				\pgfpatharc{0}{70}{0.14\myarrowsize}
        \pgfpatharc{-80}{-169.5}{4\myarrowsize}
        \pgfpatharc{150}{189}{0.95\myarrowsize and 0.95\myarrowsize}
        \pgfpatharc{0}{-40}{15\myarrowsize}
        \pgfpathmoveto{\pgfpoint{0.43\myarrowsize}{0}}
        \pgfpatharc{0}{-70}{0.14\myarrowsize}
        \pgfpatharc{80}{169.5}{4\myarrowsize}
        \pgfpatharc{-150}{-189}{0.95\myarrowsize and 0.95\myarrowsize}
        \pgfpatharc{0}{40}{15\myarrowsize}
			  \pgfpathmoveto{\pgfpoint{-7\myarrowsize}{0}}
        \pgfpatharc{0}{-70}{0.14\myarrowsize}
        \pgfpatharc{80}{169.5}{4\myarrowsize}
        \pgfpatharc{-150}{-189}{0.95\myarrowsize and 0.95\myarrowsize}
        \pgfpatharc{0}{40}{15\myarrowsize}
        \pgfpathclose
        \pgfsetstrokeopacity{0.25}
        \pgfusepathqfillstroke
        \pgfusepathqfillstroke
    }

    \pgfarrowsdeclare{myhook}{myhook}{
        \setlength{\myarrowsize}{0.6pt}
        \addtolength{\myarrowsize}{.5\pgflinewidth}
        \pgfarrowsleftextend{-4\myarrowsize-.5\pgflinewidth} 
        \pgfarrowsrightextend{.7\pgflinewidth}
    }{
        \setlength{\myarrowsize}{0.6pt} 
        \addtolength{\myarrowsize}{.5\pgflinewidth}  
        \pgfsetdash{}{+0pt}
        \pgfsetroundcap
        \pgfpathmoveto{\pgfqpoint{-2pt}{-6\pgflinewidth}}
        \pgfpathcurveto
            {\pgfqpoint{4\pgflinewidth}{-4.667\pgflinewidth}}
            {\pgfqpoint{4\pgflinewidth}{0pt}}
            {\pgfpointorigin}
        \pgfusepathqstroke
    }

		\pgfarrowsdeclare{my to}{my to}
{
  \pgfarrowsleftextend{-2\pgflinewidth}
  \pgfarrowsrightextend{\pgflinewidth}
}
{
  \pgfsetlinewidth{0.8\pgflinewidth}
  \pgfsetdash{}{0pt}
  \pgfsetroundcap
  \pgfsetroundjoin
  \pgfpathmoveto{\pgfpoint{-5.5\pgflinewidth}{7.5\pgflinewidth}}
  \pgfpathcurveto
  {\pgfpoint{-4.0\pgflinewidth}{0.1\pgflinewidth}}
  {\pgfpoint{0pt}{0.25\pgflinewidth}}
  {\pgfpoint{0.75\pgflinewidth}{0pt}}
  \pgfpathcurveto
  {\pgfpoint{0pt}{-0.25\pgflinewidth}}
  {\pgfpoint{-4.0\pgflinewidth}{-0.1\pgflinewidth}}
  {\pgfpoint{-5.5\pgflinewidth}{-7.5\pgflinewidth}}
  \pgfusepathqstroke
}

\tikzstyle{vecArrow} = [thick, decoration={markings,mark=at position
   1 with {\arrow[semithick]{open triangle 60}}},
   double distance=1.4pt, shorten >= 5.5pt,
   preaction = {decorate},
   postaction = {draw,line width=1.4pt, white,shorten >= 4.5pt}]
\tikzstyle{innerWhite} = [semithick, white,line width=1.4pt, shorten >= 4.5pt]

	\makeatletter
	\newcommand{\reqnomode}{\tagsleft@false}
	\newcommand\POSITION[3]{%
	\begingroup
	\@tempdim@x=0cm
	\@tempdim@y=\paperheight
	\advance\@tempdim@x#1
	\advance\@tempdim@y-#2
	\put(\LenToUnit{\@tempdim@x},\LenToUnit{\@tempdim@y}){#3}%
	\endgroup
	}
\makeatother

\addto\captionsenglish{}

\begin{document}
	\maketitle
	
	\begin{abstract}
		We define a tower of affine Temperley-Lieb algebras of type $\tilde{A}$. We prove that there exists a unique Markov trace on this tower, this trace comes from the Markov-Ocneanu-Jones trace on the tower of Temperley-Lieb algebras of type $A_{n}$. We define an invariant of special kind of links as an application of this trace.
	\end{abstract}

	
	
	
	\section{Introduction}
	
	About 30 years ago, V. Jones discovered one of the most famous invariants of oriented knots and links \cite{Jones_1985}. The remarkable feature of the construction of this invariant is that it works in a purely algebraic setting: it arises from certain trace functions (Markov traces) on Temperley-Lieb algebras of type $A$. Later V. Jones, himself, redefined his traces on Iwahori-Hecke algebras  of type $A$, with a parameter "z" in the ground ring.  Temperley-Lieb algebras of type $A$ are quotients of  Iwahori-Hecke algebras  of type $A$, which themselves are finite-dimensional quotients of the group algebras of Artin's braid groups (braid groups of type $A$). The other remarkable feature of the above-mentioned invariant is that it is easy to compute, yet it  does not "distinguish" totally between  oriented links, in other words there exist at least two nonequivalent links having the same Jones polynomial.  \\
	
	Since then several generalizations of the construction of Markov traces beyond the type $A$ have been achieved, yet have been restricted to the "spherical cases", that is, braid groups, Iwahori-Hecke algebras and  Temperley-Lieb algebras associated with finite Weyl groups,  see   \cite{tomDieck} for a related study in   type $B$.  A classification of Markov traces on Iwahori-Hecke algebras of type $B$ and $D$ was given by Geck and Lambropoulou in  \cite{Geck_Lambropoulou_1997} from which we can determine the Markov traces factoring by the corresponding Temperley-Lieb algebras by noticing that the space of all trace functions 
	on a finite-dimensional Iwahori-Hecke algebra is spanned by the characters of the irreducible representations of this algebra, these characters are essentially known  \cite{Geck_Pfeiffer_2000}.   \\
	
	In the affine  cases one faces new difficulties, notably   the fact  that one must deal with infinite dimensional algebras. 
	The $\tilde{A}$-type affine braid group is the braid group under question in this work; we call its elements affine braids. Geometrically, one can see several presentations in the literature, among which we choose the one corresponding to  the  $B$-type braid group (we call its elements $B$-braids)\cite{Sadek_2015_2}. In  \cite{Graham_Lehrer_2003}  we see that the $B$-type braid group is a semi-direct product of the $\tilde{A}$-type affine braid group with a normal subgroup generated by  one element (acting on the affine braids as the "Dynkin automorphism"  which is to be defined in the third section). In particular  every affine braid is a $B$-braid.  
	
	In  \cite{Sadek_2015_2}   we give a definition of an affine oriented link to be the closure of an affine braid seen as a $B$-braid;  
	
	or equivalently to be:  an oriented link in a solid torus which makes the same number of positive and negative rounds around the "middle hole"; 
	
	or equivalently to be:  an oriented braid in $S^3$  in which the number of positive rounds equals the number of negative rounds, around a fixed circle  (links which do not make rounds are counted here, in other terms every oriented link -in the usual sense- is an affine oriented link). 
	
	Hence, an invariant of oriented affine links is an invariant of oriented links.  \\  
	
	In this paper, we will essentially work on the   images of the affine braids in the affine Temperley-Lieb algebra of type  $\tilde{A}$ (in the literature  it is sometimes called non-extended affine Temperley-Lieb algebra in order to distinguish it from the extended  affine Temperley-Lieb algebra which is a slightly larger structure, see   \cite{Graham_Lehrer_2003} and   \cite{Graham_Lehrer_1998}). We will  use purely algebraic tools to establish the existence and unicity of an affine Markov trace 
	on the tower of affine Temperley-Lieb algebras.  
	
	This affine Markov trace thus 
	 defines the unique Jones-like invariant of oriented affine links, when composed with the following path:\\

\centerline{ Affine oriented links $    \longrightarrow $ Affine braid groups $   \longrightarrow   $ Affine T-L algebras.   } 

\bigskip
 
	The paper is organized as follows. 
	
	In section 2 we   recall the definition  of the affine Temperley-Lieb algebra $\widehat{TL}_{n+1} (q)$, of type 
	$\tilde A_n$. This algebra  has a basis $(g_w)_{w \in W^c(\tilde A_n)}$ indexed by the fully commutative elements 
	in the affine Coxeter group $W(\tilde A_n)$ of type 	$\tilde A_n$. 
	We also state  Jones's theorem on existence and unicity of the Markov trace on the tower of Temperley-Lieb algebras of type $A$. 
	
	In section 3 we
  build  a tower of  affine Temperley-Lieb algebras  of type  $\tilde{A}$: 
$$
			\widehat{TL}_{1}(q) \stackrel{F_{1}}{\longrightarrow}  \widehat{TL}_{2}(q) \stackrel{F_{2}} {\longrightarrow}\widehat{TL}_{3}(q) \longrightarrow ~~...~~ \widehat{TL}_{n}(q) \stackrel{F_{n}} {\longrightarrow}\widehat{TL}_{n+1}(q) \longrightarrow ...  
$$
(Remark that   $W(\tilde A_n)$   is not a parabolic subgroup of 	$W(\tilde A_{n+1})$ hence the definition of the morphism $F_n$ is not straightforward. Also remark that we do not know whether $F_n$ is injective.)  
We then  define what should be the Markov conditions in the affine case. 	
  
  We start section 4 with a study of the traces on $\widehat{TL}_{2} (q)$ and 
  $\widehat{TL}_{3} (q)$ that are invariant under the Dynkin automorphism. We  proceed and define Markov elements: those elements of $\widehat{TL}_{n+1} (q)$
  that belong either to $F_n(\widehat{TL}_{n}(q))$ or to    
  $F_n(\widehat{TL}_{n}(q)) g_{\sigma_n}  F_n(\widehat{TL}_{n}(q))$, 
  where  $\{ \sigma_1, \dots, \sigma_n\}$ is the set of Coxeter generators of 
  $W( A_n)$.    
	We then use the main result of \cite{Sadek_2013_2}, which is a classification of  fully commutative elements in the affine Weyl groups of type $\tilde{A}$,  to prove 
	that any trace on $\widehat{TL}_{n+1} (q)$ for $n\ge 2$ is uniquely defined by its values on Markov elements. This is Theorem \ref{5_1_1}, a crucial step in the paper. 
	 
	 The tower of affine Temperley-Lieb algebras of type  $\tilde{A}$  surjects onto 
the tower of Temperley-Lieb algebras of type $A$ and we prove in the beginning of section 5 that composing the Markov trace on the latter with this surjection provides indeed an affine Markov trace on the former. We then show the uniqueness of the components 
of an affine Markov trace on 	 $\widehat{TL}_{2} (q)$ and 
  $\widehat{TL}_{3} (q)$, which, combined with Theorem \ref{5_1_1}, leads us 
  to the main Theorem \ref{5_2_8} asserting 
 the existence and uniqueness of the Markov trace  on the tower of affine Temperley-Lieb algebras of type  $\tilde{A}$.

	\section{The affine Temperley-Lieb algebra}
		Let $K$ be an integral domain of characteristic $0$. Suppose that $q$ is a square invertible element in $K$ of which we fix a root  $\sqrt{q}$. For $x,y$ in a given ring we define $V(x,y):= xyx+xy+yx+x+y+1$. We mean by algebra in what follows $K$-algebra.

		We denote by $B(\tilde{A_{n}})$ (resp. $W(\tilde{A_{n}})$) the affine braid (resp. affine Coxeter) group with $n+1$ generators of type $\tilde{A}$, while we denote by $B(A_{n})$ (resp. $W(A_{n})$) the braid (resp. Coxeter) group with $n$ generators of type $A$, where $n \geq 0 $. 
		By definition $B(\tilde{A_{n}}) $ has $ \left\{ \sigma_{1}, \sigma_{2}, \dots,   \sigma_{n}, a_{n+1}   \right\}$ as a set of generators together with the following defining relations:  \\
		        	
			\begin{itemize}[label=$\bullet$, font=\normalsize, font=\color{black}, leftmargin=2cm, parsep=0cm, itemsep=0.25cm, topsep=0cm]
				 \item[(1)] $\sigma_{i} \sigma_{j} =\sigma_{j} \sigma_{i} $  where $1\leq i,j\leq n$ when $ \left| i-j\right| \geq 2$,
				 \item[(2)] $\sigma_{i}\sigma_{i+1}\sigma_{i} = \sigma_{i+1}\sigma_{i}\sigma_{i+1}$ when $1\leq i\leq n-1$,
				 \item[(3)] $\sigma_{i} a_{n+1} = a_{n+1} \sigma_{i} $ when $2\leq i \leq n-1$, 
				 \item[(4)] $\sigma_{1} a_{n+1} \sigma_{1} = a_{n+1} \sigma_{1} a_{n+1}$ for $n\geq 2$,
				 \item[(5)] $\sigma_{n} a_{n+1} \sigma_{n} = a_{n+1} \sigma_{n} a_{n+1}$ for $n\geq 2$,\\
			\end{itemize}
while $B(A_{n})$ is   generated by $ \left\{ \sigma_{1}, \sigma_{2}, \dots,  \sigma_{n}  \right\}$. 
We let $W^{c}(\tilde{A_{n}})$ (resp. $W^{c}(A_{n}))$ be the set of fully commutative elements in $W(\tilde{A_{n}})$ (resp. $W(A_{n}))$.\\

\reqnomode
		
		For $ n\geq 2$, we define $\widehat{TL}_{n+1} (q)$ to be the algebra with unit given by a set of generators  $\left\{g_{\sigma_{1}}, ...,~ g_{\sigma_{n}}, g_{a_{n+1}}\right\}$, with the following relations \cite{Graham_Lehrer_1998}:\\
	\begin{equation*}\label{definingrelations} 	
	\text{(2.0)}  \quad \left\{ \quad  \begin{aligned}
		 &g_{\sigma_{i}} g_{\sigma_{j}} =g_{\sigma_{j}} g_{\sigma_{i}}  \text{ for } 1\leq i,j\leq n \text{ and } \left| i-j\right| \geq 2, \\
		 &g_{\sigma_{i}} g_{a_{n+1}} =g_{a_{n+1}} g_{\sigma_{i}}   \text{ for  }  2\leq i \leq n-1, \\
			  &g_{\sigma_{i}}g_{\sigma_{i+1}}g_{\sigma_{i}} = g_{\sigma_{i+1}}g_{\sigma_{i}}g_{\sigma_{i+1}}  \text{ for }  1\leq i\leq n-1,\\
			   &g_{\sigma_{i}}g_{a_{n+1}}g_{\sigma_{i}} = g_{a_{n+1}}g_{\sigma_{i}}g_{a_{n+1}}  \text{ for } i= 1, n , \\
			  &g^{2}_{\sigma_{i}} = (q-1)g_{\sigma_{i}} +q  \text{ for }  1\leq i\leq n,\\
			& g^{2}_{a_{n+1}} = (q-1)g_{a_{n+1}} +q , \\ 
			  &V(g_{\sigma_{i}},g_{\sigma_{i+1}})=V(g_{\sigma_{1}},g_{a_{n+1}}) = V(g_{\sigma_{n}},g_{a_{n+1}})= 0  \text{  for }1\leq i\leq n-1.  
			  \end{aligned} \right.     \qquad 
		\end{equation*}
		
		\smallskip
We set $\widehat{TL}_{1} (q) = K$. For $n=1$, the 		
			  algebra $\widehat{TL}_{2}(q) $ is generated by two elements: $ g_{\sigma_{1}}, g_{a_{2}}$, with only  Hecke quadratic relations. That is:\\
			\begin{eqnarray}
	 			g_{\sigma_{1}}^{2} &=& (q-1) g_{\sigma_{1}}  +q \quad  \text{ and }  
				\quad g_{a_{2}}^{2} = (q-1) g_{a_{2}} +q.  \nonumber\\\nonumber
			\end{eqnarray}

		\smallskip
		 The set $\left\{ g_{w}: w \in W^{c}(\tilde{A_{n}})\right\}$ is well defined in the usual sense of the theory of Hecke algebra and it is a $K$-basis \cite[\S 2]{Fan_1996}. We set $T_{a_{n+1}}$ (resp. $T_{\sigma_{i}}$ for $1 \leq i \leq n$) to be $\sqrt{q}g_{a_{n+1}}$ (resp. $\sqrt{q}g_{\sigma_{i}}$ for $1 \leq i \leq n$). Hence, $T_{w}$ is well defined for $w \in W^{c}(\tilde{A_{n}})$, it equals $q^{\frac{l(w)}{2}}g_{w}$. The multiplication associated to the basis  $\left\{ T_{w}: w \in W^{c}(\tilde{A_{n}})\right\}$, is given as follows:
		 
		\begin{eqnarray}
			T_{w} T_{v} &=& T_{wv} ~~~~~~~~~~~~~~~~~~~~~~~~~~~~~~ \text{whenever } l(wv) = l(w) + l (v),\nonumber\\
			T_{s} T_{w} &=& \sqrt{q}(q-1) T_{w} +q^{2}  T_{sw} ~~~~~~~~ \text{whenever } l(sw) = l(w) - 1, \nonumber
		\end{eqnarray}
 		
		for   
		$w,v, wv$  in $W^{c}(\tilde{A_{n}})$ and $s$ in $\left\{\sigma_{1}, \dots, \sigma_{n}, a_{n+1}\right\}$.\\

		 In what follows we suppose that $q+1$ is  invertible in $K$, we set $\delta = \frac{1}{2+q+q^{-1}} = \frac{q}{(1+q)^{2}}$ in $K$. In view of \cite{Graham_Lehrer_2003} , for $1 \leq i \leq n$ we set $ f_{\sigma_{i}}:= \frac{ g_{\sigma_{i}}+1}{q+1}$ and $ f_{a_{n+1}}:= \frac{g_{a_{n+1}}+1}{q+1}$. In other terms $ g_{\sigma_{i}} = (q+1) f_{\sigma_{i}} -1$, and $ g_{a_{n+1}} = (q+1) f_{a_{n+1}} -1$. The set $\left\{ f_{w}: w \in W^{c}(\tilde{A_{n}})\right\}$ is well defined and it is a $K$-basis for $\widehat{TL}_{n+1} (q)$.\\
		
The classical Temperley-Lieb algebra of type $A$ with $n$ generators, $TL_{n}(q)$, can be regarded as the subalgebra of $\widehat{TL}_{n+1} (q)$ generated by  $\left\{g_{\sigma_{1}} ,...,~ g_{\sigma_{n}}\right\}$, with $\left\{ g_{w}: w \in W^{c}(A_{n})\right\}$ as $K$-basis. We set $TL_{0}(q) = K$. \\
		
		Now  we consider the following tower: 
		
		\begin{eqnarray}
				TL_{0}(q)~\subset TL_{1}(q) ~~...\subset TL_{n-1}(q) ~~\subset TL_{n}(q) ~~... \nonumber\\\nonumber
			\end{eqnarray}
			
			A fundamental result of Jones is the following: 
			
		\begin{theoreme} \label{1_1}\cite{Jones_1985}	
		There is a unique collection of traces $(\tau_{n+1})_{0 \leq n}$ on $(TL_{n})_{0\leq n}$, such that: \\
			\begin{itemize}[label=$\bullet$, font=\normalsize, font=\color{black}, leftmargin=2cm, parsep=0cm, itemsep=0.25cm, topsep=0cm]
				\item[$(1)$] $\tau_{1}(1) = 1 $.
				\item[$(2)$] For $ 1 \leq n$, we have $\tau_{n+1}(hT^{\pm1}_{\sigma_{n}}) = \tau_{n}(h) $, for any $h$ in $TL_{n-1}(q)$.
			\end{itemize}	
		\end{theoreme}		
		The collection  $(\tau_{n+1})_{0 \leq n}$ is called a Markov trace. Moreover, for any $a,b$ and $c$ in $TL_{n}(q)$ and for $ n \geq 1$, every $\tau_{n+1}: TL_{n}(q) \longrightarrow K$ verifies: 
					\begin{eqnarray}
						\tau_{n+1}(bT_{\sigma_{n}}c)= \tau_{n}(bc)$ and $\tau_{n+1}(a)=- \frac{1+q}{\sqrt{q}}\tau_{n}(a). \nonumber
					\end{eqnarray}


	\section{The tower of affine Temperley-Lieb algebras and affine Markov trace}
	
	In this section we define a tower of affine Temperley-Lieb algebras, we show that this tower "surjects" onto the tower of Temperley-Lieb algebras mentioned in the previous section, and we define the affine Markov trace.\\ 
	
	We consider the Dynkin diagram of the group $B(\tilde{A_{n}})$. 			
					
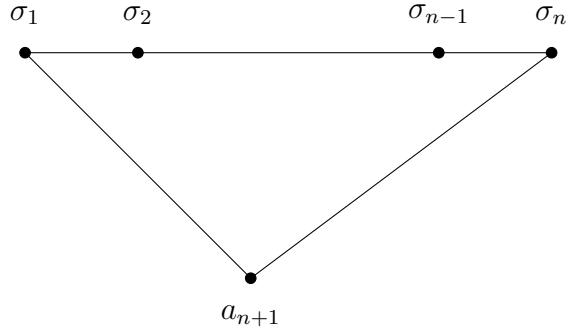
			\begin{figure}[ht]
				\centering
				\begin{tikzpicture}

 \node at (0,0.5) {$\sigma_{1}$}; 
  \filldraw (0,0) circle (2pt);
   
  \draw (0,0) -- (1.5, 0);
  
  \node at (1.5,0.5) {$\sigma_{2}$};
  \filldraw (1.5,0) circle (2pt);

  \draw (1.5,0) -- (5.5, 0);

  \node at (5.5,0.5) {$\sigma_{n-1}$};
  \filldraw (5.5,0) circle (2pt);
 
  \draw (5.5,0) -- (7, 0);
  
  \node at (7,0.5) {$\sigma_{n}$};
  \filldraw (7,0) circle (2pt);

  \draw (7,0) -- (3, -3);
  
  \filldraw (3, -3) circle (2pt);
  \node at (3, -3.5) {$a_{n+1}$};

  \draw (3, -3) -- (0, 0);
               \end{tikzpicture}
			\end{figure}
	
We denote the Dynkin automorphism $(\sigma_{1}\mapsto \sigma_{2} \mapsto \cdots \mapsto  \sigma_{n} \mapsto a_{n+1} \mapsto \sigma_{1})$ by  $\psi_{n+1}$.	We have the following injection:
	
	\begin{eqnarray}
					G_{n}: K[B(\tilde{A_{n-1}})] &\longrightarrow& K[B(\tilde{A_{n}})]  \nonumber\\
					\sigma_{i} &\longmapsto& \sigma_{i}$ ~~~ \text{for} $1\leq i\leq n-1 \nonumber\\
					a_{n} &\longmapsto& \sigma_{n} a_{n+1}\sigma^{-1}_{n} \nonumber
				\end{eqnarray}

	We prove in \cite{Sadek_Thesis}, to which we refer for details, that 
	$G_n$ induces an injective homomorphism of   the corresponding Hecke algebras \cite[Proposition 4.3.3]{Sadek_Thesis}. This homomorphism induces in turn a map at the Temperley-Lieb level that is an algebra homomorphism \cite[\S 5.2.3]{Sadek_Thesis},  but  in the affine case,  the possible lack of injectivity forces us to use different notations for the generators of $\widehat{TL}_{n}(q)$ and $\widehat{TL}_{n+1}(q)$. In the following Proposition we use 
	  $\left\{t_{\sigma_{1}}, ..., t_{\sigma_{n-1}}, t_{a_{n}}\right\}$
as the	set of generators for $\widehat{TL}_{n}(q)$ satisfying relations (2.0) (with  $t $ replacing   $g $   and $n$ replacing $n+1$).   
 \\

	
	\begin{proposition} \label{2_1} The injection $G_{n}$ induces the following morphism of algebras:
	
			\begin{eqnarray}
				F_{n}: \widehat{TL}_{n}(q) &\longrightarrow& \widehat{TL}_{n+1}(q) \nonumber\\
				t_{\sigma_{i}} &\longmapsto & g_{\sigma_{i}} \text{ for } 1 \leq i \leq n-1 \nonumber\\
				t_{a_{n}} &\longmapsto & g_{\sigma_{n}} g_{a_{n+1}} g^{-1}_{\sigma_{n}}. \nonumber
			\end{eqnarray}
	\end{proposition}	 
	
		 Notice that the conjugation by $ \sigma_{n} \sigma_{n-1}.. \sigma_{1}a_{n+1} $ in  $ B(\tilde{A_{n}}) $  acts   on the image of $ B(\tilde{A_{n-1}}) $ as 
		 $\psi_{n}^{-1}$:  ($\sigma_{1} \mapsto a_{n} \mapsto \sigma_{n-1} \mapsto \sigma_{n-2} \mapsto.. \sigma_{2} \mapsto \sigma_{1}$). We denote by $\Psi$ the automorphism of 
				$ G_{n}(B(\tilde{A_{n-1}}) )$ image of $\psi_n$. 			 				We thus 
				 write $(\sigma_{n} .. \sigma_{1} a_{n+1})^{d}h = \Psi^{-d} (h)  (\sigma_{n} .. \sigma_{1} a_{n+1})^{d} $, for any $h$ in $ G_{n}(B(\tilde{A_{n-1}}) )$. 
				 This convention is specially important when used   for the affine Temperley-Lieb algebra.	Indeed   $\psi_{n+1}$ also acts as an automorphism of the algebra  $ \widehat{TL}_{n+1}(q) $ and the conjugation by 
				 $g_{ \sigma_{n} \sigma_{n-1}.. \sigma_{1}a_{n+1} }$  
			 in $ \widehat{TL}_{n+1}(q) $ stabilizes the image 
			 of $\widehat{TL}_{n}(q)$, on which it acts like the image 
			 $\Psi$ of $\psi_{n}$ under $F_n$. We then write likewise:

			 	\begin{equation}\label{conjugacy}
				g_{ \sigma_{n} \sigma_{n-1}.. \sigma_{1}a_{n+1} }^d \,  h = 
			  \Psi^{-d} (h)   \, 
			g_{ \sigma_{n} \sigma_{n-1}.. \sigma_{1}a_{n+1} }^d  
		\quad	\text{ for all } h \in  F_n(\widehat{TL}_{n}(q)) 
		\text{ and } d \in \mathbb  Z .  
		\end{equation}
		\\

		We also prove   in \cite[\S 5.2.2]{Sadek_Thesis}: \\
		
	\begin{proposition}\label{En}   The following map is a surjection of algebras
	
			\begin{eqnarray}
				E_{n}: \widehat{TL}_{n+1}(q) &\longrightarrow& TL_{n}(q) \nonumber\\
				g_{\sigma_{i}} &\longmapsto & g_{\sigma_{i}} \text{ for } 1 \leq i \leq n \nonumber\\
				g_{a_{n+1}} &\longmapsto & g_{\sigma_{1}} ... g_{\sigma_{n-1}} g_{\sigma_{n}} g^{-1}_{\sigma_{n-1}} ... g^{-1}_{\sigma_{1}}. \nonumber
			\end{eqnarray}
	
	Moreover, the following diagram commutes: 
	
	\begin{center}    

	\begin{tikzpicture}

			\matrix[matrix of math nodes,row sep=1cm,column sep=1cm]{
			|(A)| \widehat{TL}_{n} (q)    & & & &    |(B)| \widehat{TL}_{n+1} (q)  \\
			                              & & & &                      \\								
			|(C)| TL_{n-1}(q)             & & & &    |(D)| TL_{n}(q)    \\
				};

				\path (A) edge[-myonto,line width=0.42pt]  node[above, xshift=-5mm, yshift=-2mm, rotate=0] {\footnotesize $E_{n-1} $}    (C);
				
				 				\path (A) edge[-myto,line width=0.42pt] node[above, xshift=1.5mm, yshift=0mm, rotate=0] {\footnotesize $F_{n}$}      (B);
				
				\path (D) edge[-myhook,line width=0.42pt] node[below, xshift=1.5mm, yshift=0mm, rotate=0] {\footnotesize $$}    (C);
				\path (C) edge[-myto,line width=0.42pt]      (D);

				\path (B) edge[-myonto,line width=0.42pt]   node[above, xshift=5mm, , yshift=-2mm, rotate=0] {\footnotesize $E_{n}$}   (D);

		\end{tikzpicture}	
	\end{center}
		\end{proposition}			

	 Note that $E_{n}$ composed with the natural inclusion of $TL_{n}(q)$ into $\widehat{TL}_{n+1}(q)$, gives $Id_{TL_{n}(q)}$. \\

			In view of Proposition \ref{2_1} we can consider the tower of affine T-L algebras (it is not known whether it is a tower of faithful arrows or not):\\ 
		\begin{eqnarray}
			\widehat{TL}_{1}(q) \stackrel{F_{1}}{\longrightarrow}  \widehat{TL}_{2}(q) \stackrel{F_{2}} {\longrightarrow}\widehat{TL}_{3}(q) \longrightarrow ~~...~~ \widehat{TL}_{n}(q) \stackrel{F_{n}} {\longrightarrow}\widehat{TL}_{n+1}(q) \longrightarrow ... \nonumber\\\nonumber
		\end{eqnarray}
				
			\begin{definition} \label{5_2_1}
				We call $(\hat{\tau}_{n})_{1 \leq n}$ an affine Markov trace, if every $\hat{\tau}_{n}$ is a trace function on $\widehat{TL}_{n} (q)$ with the following conditions:\\
			
				\begin{itemize}
					\item $\hat{\tau}_{1}(1) = 1$, 
										\item $\hat{\tau}_{n+1}(F_{n}(h)T^{\pm1}_{\sigma_{n}}) =  \hat{\tau}_{n}(h)$, for all $h \in  \widehat{TL}_{n} (q)$ and for $n \geq 1$.
					\item $\hat{\tau}_{n}$ is invariant under the Dynkin automorphism $\psi_{n}$ for all $n$.\\               
				\end{itemize}
			\end{definition}
		
		\begin{remarque} 
				We notice that the second condition gives us that $\hat{\tau}_{n+1}\big(F_{n}(h)T^{-1}_{\sigma_{n}}\big) =  \hat{\tau}_{n}\big(h\big)$, which means that:
				\begin{eqnarray}
					\hat{\tau}_{n+1}\big(F_{n}(h)[\frac{1}{q^{2}}T_{\sigma_{n}}- \frac{q-1}{q\sqrt{q}}] \big) =  \hat{\tau}_{n}\big(h). \text{ Thus } \hat{\tau}_{n+1} \big(F_{n}(h)\big) = -\frac{q+1}{\sqrt{q}} \hat{\tau}_{n}\big(h\big). \nonumber\\\nonumber
				\end{eqnarray}
			\end{remarque}
			
			\begin{remarque}
			
			The third condition of Definition \ref{5_2_1} is, in fact, not independent, i.e., it results from the first and second conditions. We just have to see that if we have two elements in $ \widehat{TL}_{n} (q)$, say $x$ and $y$, such that $\psi_{n} (x) = y$, then $F_{n}(x)$ and $ F_{n}(y)$ are conjugate in $ \widehat{TL}_{n+1} (q)$, by some power of the element $g_{\sigma_{n} ...~\sigma_{1}a_{n+1}} $, which results from (1) above. Nevertheless, we will keep viewing it as a condition and study traces invariant under the Dynkin automorphism.  \\
	
	      \end{remarque}
		

	
	\section{On the space of traces on $\widehat{TL}_{n+1}(q)$} 
	
	\bigskip
	
	\subsection{Traces on $\widehat{TL}_{2} (q)$}

			The algebra $\widehat{TL}_{2}(q) $ is generated by the two elements  $ {f}_{\sigma_{1}}$ and $ {f}_{a_{2}}$ with relations ${f}^{2}_{\sigma_{1}} =  {f}_{\sigma_{1}}$ and ${f}^{2}_{a_{2}} =  {f}_{a_{2}}$. Moreover, $\widehat{TL}_{2}(q)$ has  $\left\{ {f}_{w}; w \in   W(\tilde{A}_{1}) \right\}$ as a $K$-basis. The aim of this subsection is to parametrize all traces over this algebra that  are invariant under the action of the Dynkin automorphism $\psi_{2}$, which exchanges   $ {f}_{\sigma_{1}}$ and $ {f}_{a_{2}}$. Clearly, any trace that has the same value on $ {f}_{\sigma_{1}}$ and $ {f}_{a_{2}}$ is invariant under the Dynkin automorphism $\psi_{2}$.\\
		    
			\begin{proposition} \label{5_1_10}
			Let  $A_{0},A_{1}$ and $(\alpha_{i})_{i\geq1} $ be arbitrary elements in the ground field. Then, there exists a unique trace $t$ on $\widehat{TL}_{2}(q)$, invariant by the action of $\psi_{2}$, in such a way that: 
			$A_{0} = t(1), ~ A_{1}= t( {f}_{\sigma_{1}})$ and $\alpha_{i}=t\big(( {f}_{\sigma_{1}a_{2}})^{i} \big)$.
			
			\end{proposition}
			
			\begin{demo}
			
			We start by the existence. Let $t$ be the linear function given by:
			\begin{eqnarray}
				t:\widehat{TL}_{2} (q) &\longrightarrow& K \nonumber\\
				t(1) &=& A_{0} \nonumber\\
				t( {f}_{\sigma_{1}}) &=& t( {f}_{a_{2}})  = A_{1} \nonumber\\
				t\big(( {f}_{\sigma_{1}a_{2}})^{s} \big) &=& t \big(( {f}_{a_{2}\sigma_{1}})^{s} \big)  =  t\big(( {f}_{\sigma_{1}a_{2}})^{s} {f}_{\sigma_{1}}\big) = t\big(( {f}_{a_{2}\sigma_{1}})^{s} {f}_{a_{2}} \big) = \alpha_{s}, \nonumber
			\end{eqnarray}
where $A_{0},A_{1}$ and $\alpha_{i} $ are arbitrary elements in the ground field for $ i \geq 1 $.\\
		 
	 		We show that this linear function is a trace. First we see that $t$ is, by definition, invariant under the Dynkin automorphism $ \psi_{2}$. In order to show that $t$ is a trace, we show that $ t(xy) = t(yx) $ for any $ x$ and $y$ in $\widehat{TL}_{2}(q)$. The way to do so, is to show that it is true when $ x $ is any element of the left column, and $y$ is any element of the right column, in the following table:\\ 

\newpage 
		     \hspace{1.1cm}$[1]( {f}_{\sigma_{1}a_{2}})^{k}$                  \hspace{7.6cm}$[1'] ( {f}_{\sigma_{1}a_{2}})^{h}$\\
				
			\hspace{1.1cm}$[2]( {f}_{a_{2}\sigma_{1}})^{k}$                  \hspace{7.6cm}$[2'] ( {f}_{a_{2}\sigma_{1}})^{h}$\\ 
		
		     \hspace{1.1cm}$[3]( {f}_{\sigma_{1}a_{2}})^{k} {f}_{\sigma_{1}}$    \hspace{7.1cm}$[3'] ( {f}_{\sigma_{1}a_{2}})^{h} {f}_{ \sigma_{1}}$\\
		     
		     \hspace{1.1cm}$[4]( {f}_{a_{2}\sigma_{1}})^{k} {f}_{a_{2}}$         \hspace{7.1cm}$[4'] ( {f}_{a_{2}\sigma_{1}})^{h} {f}_{a_{2}}$\\
		
		     The only cases to consider are \textbf{[1-2'], [1-3'], [1-4'] and  [3-4']}, up to applying $\psi_{2}$.\\

			\textbf{[1-2']}:			
\begin{eqnarray}
\text{Here, }t\big(xy\big) &=&  t\big(( {f}_{\sigma_{1}a_{2}})^{k} ( {f}_{a_{2}\sigma_{1}})^{h}\big) =  t\big(( {f}_{\sigma_{1}a_{2}})^{k} {f}_{a_{2}\sigma_{1}} ( {f}_{a_{2}\sigma_{1}})^{h-1}\big) = t\big(( {f}_{\sigma_{1}a_{2}})^{k} ( {f}_{\sigma_{1}a_{2}})^{h-1} {f}_{\sigma_{1}}\big)  \nonumber\\\nonumber\\
				&=& t\big(( {f}_{\sigma_{1}a_{2}})^{k+h-1} {f}_{\sigma_{1}}\big) = \alpha_{k+h-1}, \nonumber
			\end{eqnarray}			
			\begin{eqnarray}
				  \text{while  }t\big(yx\big) &=&  t\big(({f}_{a_{2}\sigma_{1}})^{h}( {f}_{\sigma_{1}a_{2}})^{k}\big)= t\big( ( {f}_{a_{2}\sigma_{1}})^{h} {f}_{\sigma_{1}a_{2}}( {f}_{\sigma_{1}a_{2}})^{k-1}\big) = t\big( ( {f}_{a_{2}\sigma_{1}})^{h}( {f}_{a_{2}\sigma_{1}})^{k-1} {f}_{a_{2}}\big)  \nonumber\\\nonumber\\
				 &=& \alpha_{k+h-1}. \nonumber
			\end{eqnarray}
	
	        	\textbf{[1-3']}:			
\begin{eqnarray}
\text{Here, }t(xy) = t\big(( {f}_{\sigma_{1}a_{2}})^{k}( {f}_{\sigma_{1}a_{2}})^{h} {f}_{ \sigma_{1}}\big) =  t\big(( {f}_{\sigma_{1}a_{2}})^{k+h} {f}_{ \sigma_{1}}\big), \text{ which is equal to } \alpha_{k+h}, \nonumber
			\end{eqnarray}
			\begin{eqnarray}
				\text{while  }t(yx) = t\big(( {f}_{\sigma_{1}a_{2}})^{h} {f}_{ \sigma_{1}}( {f}_{\sigma_{1}a_{2}})^{k} \big) = t\big(( {f}_{\sigma_{1}a_{2}})^{h+k}\big) = \alpha_{k+h}. \nonumber
			\end{eqnarray}

			\textbf{[1-4']}:			
\begin{eqnarray}
\text{Here, }t\big(xy\big)&=&t\big(( {f}_{\sigma_{1}a_{2}})^{k} ( {f}_{a_{2}\sigma_{1}})^{h} {f}_{a_{2}}\big) = t\big(( {f}_{\sigma_{1}a_{2}})^{k} {f}_{a_{2}} ( {f}_{\sigma_{1}a_{2}})^{h}\big) = t\big(( {f}_{\sigma_{1}a_{2}})^{k+h}\big) =\alpha_{k+h}, \nonumber\\\nonumber\\
				\text{while  }t\big(yx\big) &=& t\big(( {f}_{a_{2}\sigma_{1}})^{h} {f}_{a_{2}}( {f}_{\sigma_{1}a_{2}})^{k} \big) = t\big( {f}_{a_{2}}( {f}_{\sigma_{1}a_{2}})^{h}( {f}_{\sigma_{1}a_{2}})^{k} \big) = t\big( {f}_{a_{2}}( {f}_{\sigma_{1}a_{2}})^{h+k}\big) = \alpha_{k+h}.\nonumber
			\end{eqnarray}
		
			\textbf{[3-4']}:
\begin{eqnarray}
				\text{We see that:  }t\big(xy\big)&=& t\big(( {f}_{\sigma_{1}a_{2}})^{k} {f}_{\sigma_{1}} ( {f}_{a_{2}\sigma_{1}})^{h} {f}_{a_{2}}\big) = t \big(( {f}_{\sigma_{1}a_{2}})^{k+h+1}\big) =\alpha_{k+h+1}, \nonumber\\\nonumber\\
				\text{while }t\big(yx\big)&=&  t\big( ( {f}_{a_{2}\sigma_{1}})^{h} {f}_{a_{2}}( {f}_{\sigma_{1}a_{2}})^{k} {f}_{\sigma_{1}}\big) =  t\big( ( {f}_{a_{2}\sigma_{1}})^{h+k+1}\big) = \alpha_{k+h+1}.\nonumber							
			\end{eqnarray}
		    
			Now, we end the proof by showing the uniqueness. Let $t$ be a $\psi_{2} $-invariant trace on $\widehat{TL}_{2} (q)$. We have necessarily $t( {f}_{\sigma_{1}}) = t( {f}_{a_{2}})$, since $t$ is a $\psi_{2} $-invariant, call this value $A_{1}$. For every $s \geq 1$ we have $t\big(( {f}_{\sigma_{1}a_{2}})^{s} \big) = t \big(( {f}_{a_{2}\sigma_{1}})^{s} \big)$, since $t$ is a trace, call this value $\alpha_{s} $. Finally, we have $\alpha_{s} =  t\big(( {f}_{\sigma_{1}a_{2}})^{s} {f}_{\sigma_{1}}\big) = t\big(( {f}_{a_{2}\sigma_{1}})^{s} {f}_{a_{2}} \big)$, since $t$ is a trace, and ${f}_{a_{2}},~ {f}_{\sigma_{1}}$ are idempotent. Call $t(1) = A_{0}$, thus, $t$ is uniquely determined by $ A_{0},~ A_{1}$ and $\alpha_{s} $, for $s \geq 1$. 			
			\end{demo}
		 \subsection{Traces on $\widehat{TL}_{3} (q)$}
		
			In this subsection, we parametrize all the traces over $\widehat{TL}_{3} (q)$, which are invariant under the action of the Dynkin automorphism $\psi_{3} $.\\
			
			The affine Temperley-Lieb algebra in three generators $g_{\sigma_{1}},g_{\sigma_{2}}$ and $g_{a_{3}}$ can be presented by those generators with the relations of Hecke algebra, together with:
			\begin{eqnarray}
				V(g_{\sigma_{1}},g_{\sigma_{2}})=V(g_{\sigma_{1}},g_{a_{3}})=V(g_{\sigma_{2}},g_{a_{3}})=0.  \nonumber	
			\end{eqnarray}
			 As we did for $\widehat{TL}_{2} (q)$, we   change the generators as in section 2: we use $f_{\sigma_{i}}= \frac{g_{\sigma_{i}}+1}{q+1}$   for $i=1,2$, the same for $f_{a_{3}}$. $\widehat{TL}_{3}(q)$ is presented by these three generators and  the following relations:
			\begin{eqnarray}
				f^{2}_{\sigma_{i}} &=& f_{\sigma_{i}} ~\text{ for } i=1,2 \text{ and } f^{2}_{a_{3}} = f_{a_{3}},  \nonumber\\\nonumber\\				
				f_{\sigma_{i}}f_{a_{3}}f_{\sigma_{i}} &=& \delta f_{\sigma_{i}} \text{ and }  f_{a_{3}}f_{\sigma_{i}}f_{a_{3}} = \delta f_{a_{3}}   ~\text{ for } i=1,2 \nonumber\\\nonumber\\
				f_{\sigma_{1}}f_{\sigma_{2}}f_{\sigma_{1}} &=& \delta f_{\sigma_{1}} \text{ and } f_{\sigma_{2}}f_{\sigma_{1}}f_{\sigma_{2}} = \delta f_{\sigma_{2}}.\nonumber\\\nonumber	
			\end{eqnarray}
Here we will use the $K$-basis $\left\{f_{w}; w \in  W^{c}(\tilde{A}_{2}) \right\}$.\\

			\begin{lemme}
				Let $h$ and $k$ be two positive integers. Then: \\
			
				\begin{equation*}
					\big(f_{\sigma_{2}\sigma_{1}a_{3}}\big)^{k} \big(f_{\sigma_{1}\sigma_{2}a_{3}}\big)^{h} =  
					\left\{ 
						\begin{array}{ccc}
							\delta^{3h}\big(f_{\sigma_{2}\sigma_{1}a_{3}}\big)^{k-h}  &\text{for  }h<k.& \\
							\\
							\delta^{3k-1}f_{\sigma_{2}a_{3}} \big( f_{\sigma_{1}\sigma_{2}a_{3}} \big)^{h-k} &\text{for  }h\geq k.&  
						\end{array}
					\right.
				\end{equation*}
				\\		
								
				\begin{equation*}
					\big(f_{\sigma_{1}\sigma_{2}a_{3}}\big)^{h}\big(f_{\sigma_{2}\sigma_{1}a_{3}}\big)^{k}= 
					\left\{ 
						\begin{array}{ccc}
							\delta^{3k}\big(f_{\sigma_{1}\sigma_{2}a_{3}}\big)^{h-k} &\text{for  } h > k.&  \\
							\\ 
							\delta^{3h-1}f_{\sigma_{1}a_{3}}\big(f_{\sigma_{2}\sigma_{1}a_{3}}\big)^{k-h}  &\text{for  }h\leq k.& \\
					\end{array}
					\right. 
		     	\end{equation*}
			\end{lemme} 
            
			\vspace{0.5 cm}
			
			\begin{demo}
			By induction, with a direct computation the Lemma follows. 
			\end{demo}
			
			\vspace{0.5 cm}
			
			Now we parametrize all the traces on $\widehat{TL}_{3}(q)$  which are invariant by the  Dynkin automorphism $\psi_{3}$. We know that any element of the $K$-basis  $\left\{f_{w}; w \in  W^{c}(\tilde{A}_{2}) \right\}$ can be written as follows (see \cite{Sadek_2013_2}):\\  
		
		%
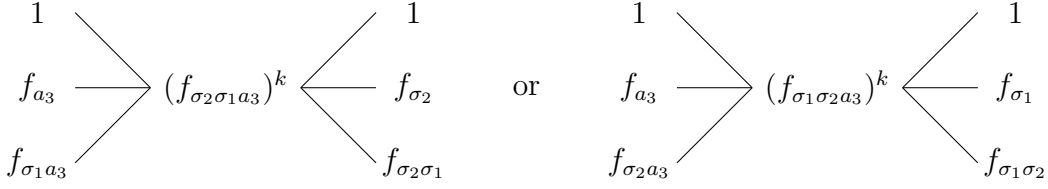
				\begin{tikzpicture}
\begin{scope}[xscale = 1]

  \node at (-1.5,1)  {$1$};
  \node at (-1.5,0)  {$f_{a_{3}}$};
  \node at (-1.5,-1) {$f_{\sigma_{1} a_{3}}$};

	\draw (-1,1)  -- (0,0);
	\draw (-1,0)  -- (0,0);
	\draw (-1,-1) -- (0,0);
	
  \node at (1,0)  {$(f_{\sigma_{2}\sigma_{1}a_{3}})^{k}$};

	\draw (2,0)  -- (3,-1);
	\draw (2,0)  -- (3,0);
	\draw (2,0)  -- (3,1);

  \node at (3.5,1)  {$1$};
  \node at (3.5,0)  {$f_{\sigma_{2}}$};
  \node at (3.5,-1) {$f_{\sigma_{2} \sigma_{1}}$};

  \node at (5,0)  {or};
\end{scope}
\begin{scope}[xscale = 1, xshift = 8cm]

  \node at (-1.5,1)  {$1$};
  \node at (-1.5,0)  {$f_{a_{3}}$};
  \node at (-1.5,-1) {$f_{\sigma_{2} a_{3}}$};

	\draw (-1,1)  -- (0,0);
	\draw (-1,0)  -- (0,0);
	\draw (-1,-1) -- (0,0);
	
  \node at (1,0)  {$(f_{\sigma_{1}\sigma_{2}a_{3}})^{k}$};

	\draw (2,0)  -- (3,-1);
	\draw (2,0)  -- (3,0);
	\draw (2,0)  -- (3,1);

  \node at (3.5,1)  {$1$};
  \node at (3.5,0)  {$f_{\sigma_{1}}$};
  \node at (3.5,-1) {$f_{\sigma_{1} \sigma_{2}}$};
\end{scope}

			
               \end{tikzpicture}

	\vspace{0.5 cm}		
			
			\begin{lemme}\label{5_1_13}
		
				Let $k$ be a positive integer, then for any $w$, such that $l(w)=3k$, the element $f_{w}$ is the image, under some power of the Dynkin automorphism $\psi_{3}$, of one of the following elements $(f_{\sigma_{2}\sigma_{1}a_{3}})^{k} $ or $(f_{\sigma_{1}\sigma_{2}a_{3}})^{k}$. Similarly for any $u$ of length $3k+1$ (resp. $3k+2$), the element $f_{u}$ is the image under a power of $\psi_{3}$ of one of the following elements $(f_{\sigma_{2}\sigma_{1}a_{3}})^{k}f_{\sigma_{2}} $ or $(f_{\sigma_{1}\sigma_{2}a_{3}})^{k}f_{\sigma_{1}}$\big(resp. $(f_{\sigma_{2}\sigma_{1}a_{3}})^{k}f_{\sigma_{2}\sigma_{1}}$ or $(f_{\sigma_{1}\sigma_{2}a_{3}})^{k}f_{\sigma_{1}\sigma_{2}}\big)$.\\
			\end{lemme}

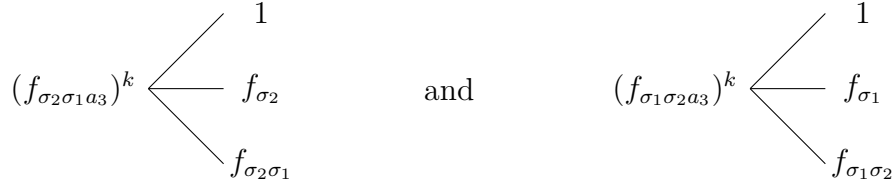
			\begin{figure}[ht]
				\centering
				\begin{tikzpicture}
               \begin{scope}[xscale = 1]
	
  \node at (1,0)  {$(f_{\sigma_{2}\sigma_{1}a_{3}})^{k}$};

	\draw (2,0)  -- (3,-1);
	\draw (2,0)  -- (3,0);
	\draw (2,0)  -- (3,1);

  \node at (3.5,1)  {$1$};
  \node at (3.5,0)  {$f_{\sigma_{2}}$};
  \node at (3.5,-1) {$f_{\sigma_{2} \sigma_{1}}$};
 
  \node at (6,0)  {and};
\end{scope}
\begin{scope}[xscale = 1, xshift = 8cm]

  \node at (1,0)  {$(f_{\sigma_{1}\sigma_{2}a_{3}})^{k}$};

	\draw (2,0)  -- (3,-1);
	\draw (2,0)  -- (3,0);
	\draw (2,0)  -- (3,1);

  \node at (3.5,1)  {$1$};
  \node at (3.5,0)  {$f_{\sigma_{1}}$};
  \node at (3.5,-1) {$f_{\sigma_{1} \sigma_{2}}$};

\end{scope}
               \end{tikzpicture}

			\end{figure}
		\begin{demo}
		The proof is direct, by induction over $k$. 
		\end{demo}
			
			 \begin{proposition}\label{3_4}
				For $i \geq 1$, let $B_{0}, B_{1},B_{2}$ and $\beta_{i}$ be in $K$. Then, there exists a unique, $\psi_{3}$-invariant, trace over $\widehat{TL}_{3}(q)$, say $s$, such that: $B_{0}=s(1),~ B_{1}=s(f_{\sigma_{1}}),~ B_{2}= s(f_{\sigma_{1}\sigma_{2}}),~ \beta_{1}= s(f_{\sigma_{1}\sigma_{2}a_{3}}),~  \beta_{k}=s((f_{\sigma_{1}\sigma_{2}a_{3}})^{k}f_{\sigma_{1}} ) $ and  $\beta_{k} = \frac{1}{\delta} s((f_{\sigma_{1}\sigma_{2}a_{3}})^{k}f_{\sigma_{1}\sigma_{2}} ) $, for $k \geq 1$.\\
			\end{proposition}

			\begin{demo}
			For the existence, we consider the following linear map $s$. We can show, using Lemma \ref{5_1_13}, that it is indeed a $\psi_{3}$-invariant trace.			
			\begin{eqnarray}
				s \text{ is given as follows, }s:\widehat{TL}_{3} (q) &\longrightarrow& K \nonumber\\
				s(1) &=& B_{0}, \nonumber\\
				s(f_{\sigma_{1}}) &=& s(f_{\sigma_{2}}) = s(f_{a_{3}})  = B_{1}, \nonumber\\
				s(f_{u}) &=& B_{2} \text{ for any } u \text{ in } W^{c}(\tilde{A}_{2}) \text{ with } l(u) = 2, \nonumber
			\end{eqnarray}
		     
			
			\begin{equation*}
				~~~~~~~~~~~~~~~~~~~~~~~~~~~~~~~~~\text{ and }s(f_{v}) = 
				\left\{ 
					\begin{array}{ccc}
							\beta_{k} ~~~~\text{when} ~~ l(v) = 3k, ~~ \text{or} &  l(v) = 3k+1, \\
							\delta \beta_{k} ~~ \text{when} ~~ l(v) = 3k+2,& \text{for}~~ k \geq 1, \\							
					\end{array}
				\right.
			\end{equation*}
	where $\beta_{k}$ (for $1 \leq k$), $B_{0}, B_{1}$ and $B_{2}$ are arbitrary in the field $K$. 
	
	For the uniqueness, we follow the steps of the proof of Proposition \ref{5_1_10}.	        
			\end{demo}
			
			\subsection{Markov elements}
		
		We consider $F_{n}:\widehat{TL}_{n} (q)\longrightarrow \widehat{TL}_{n+1} (q)$ of Proposition \ref{2_1}. In this subsection we set $F:= F_{n}$. We give a definition of Markov elements in $\widehat{TL}_{n+1} (q)$ for $ 2 \leq n $. Then we show that any trace over $\widehat{TL}_{n+1} (q)$ is uniquely determined by its values on those elements.

			\begin{definition}
				For $F$ as above, and $n \geq 2$, a Markov element in $\widehat{TL}_{n+1} (q)$ is any element of the form $ A g^{\epsilon}_{\sigma_{n}} B $, where $A$ and $B$ are in $F(\widehat{TL}_{n} (q))$ and $\epsilon \in \left\{ 0,1 \right\}$. 
			\end{definition}
		 
		The aim of this subsection is to prove the following Theorem.

			\vspace{0.25cm}		
			
			\begin{theoreme} \label{5_1_1} 
				Let $\tau_{n+1}$ be any trace over $\widehat{TL}_{n+1} (q)$ for $n \geq 2$. Then, $\tau_{n+1}$ is uniquely defined by its values on the Markov elements in $\widehat{TL}_{n+1} (q)$.\\		
             \end{theoreme}

			The proof of Theorem \ref{5_1_1} is divided into two parts. In the first we show some general facts, in the second we prove the above Theorem  for $ n \geq 3  $. For $ n=2$ we will not give the proof, as it  is pretty long and 
	  available in \cite{Sadek_2015_1}. Note that in the final proof of 
	  Theorem \ref{5_2_8} we will not need the case $n = 2$ above.  \\
			
		 We use in the proof the same notational convention as in Proposition \ref{2_1}:  we denote by 
	  $\left\{t_{\sigma_{1}}, ..., t_{\sigma_{n-1}}, t_{a_{n}}\right\}$
  the	set of generators for $\widehat{TL}_{n}(q)$ satisfying relations (2.0) (with  $t $ replacing   $g $   and $n$ replacing $n+1$), and we use the corresponding basis 
		$\left\{ t_{w}: w \in W^{c}(\tilde{A_{n-1}})\right\}$. \\

			\textbf{Part 1 }\\
	
	In this part, we suppose that $\tau_{n+1}$ is any trace on $\widehat{TL}_{n+1} (q)$. We will apply $\tau_{n+1}$ to $\widehat{TL}_{n+1} (q)$ assuming that $2 \leq n$, and show that  $\tau_{n+1}$ is uniquely determined on $\widehat{TL}_{n+1} (q)$ by its values on the positive powers of $g_{\sigma_{n} \sigma_{n-1} ..\sigma_{1} a_{n+1}}$, in addition to its values on Markov elements. From now on we denote by $w$: an arbitrary element in $W^{c}(\tilde{A_{n}})$.
	
			\begin{lemme}\label{5_1_3}
				In $\widehat{TL}_{n+1} (q) $ we have: \\

					\vspace{-0.8cm}
					\begin{eqnarray}
						(1)~g_{\sigma_{n}} (g_{\sigma_{n} \sigma_{n-1} ..\sigma_{1} a_{n+1}})^{k} &=& (q-1) (g_{\sigma_{n} \sigma_{n-1} ..\sigma_{1} a_{n+1}})^{k} + \sum\limits^{i=k-1}_{i=1} f_{i}  (g_{\sigma_{n} \sigma_{n-1} ..\sigma_{1} a_{n+1}})^{i} \nonumber \\
						& & + A \big(g_{\sigma_{n-1} \sigma_{n-2} ..\sigma_{1}}F(t_{a_{n}})\big)^{k} g_{\sigma_{n}}\prod^{j=k-1}_{j=0} \Psi^{j} \big(F((t_{a_{n}})^{-1})\big), \nonumber
					\end{eqnarray}

					\vspace{-0.5cm}
					\begin{eqnarray}
						(2)~(g_{\sigma_{n} \sigma_{n-1} ..\sigma_{1} a_{n+1}})^{k}g_{\sigma_{n}} &=& (q-1) (g_{\sigma_{n} \sigma_{n-1} ..\sigma_{1} a_{n+1}})^{k} + \sum^{i=k-1}_{i=1} h_{i}  (g_{\sigma_{n} \sigma_{n-1} ..\sigma_{1} a_{n+1}})^{i}\nonumber \\
						& & + A \prod^{j=k-1}_{j=0} \Psi^{-j} \big((g_{\sigma_{{n-1}}})^{-1}\big) g_{\sigma_{n}}\big(g_{\sigma_{n-1} \sigma_{n-2} ..\sigma_{1}}F(t_{a_{n}})\big)^{k}, \nonumber
					\end{eqnarray}
	
				with $ A $ in the ground field and  $ f_{i},h_{i} $ in $F (\widehat{TL}_{n} (q)) $.
					
			\end{lemme}
			
			\begin{demo}
				\begin{eqnarray}
					g_{\sigma_{n}} \big(g_{\sigma_{n} \sigma_{n-1} ..\sigma_{1} a_{n+1}}\big)^{k} &=& \big(q-1\big) \big(g_{\sigma_{n} \sigma_{n-1} ..\sigma_{1} a_{n+1}}\big)^{k} \nonumber\\
					& & + q g_{\sigma_{n-1} \sigma_{n-2} ..\sigma_{1}}F\big(t_{a_{n}}\big) g_{\sigma_{n}} F\big((t_{a_{n}})^{-1}\big) \big(g_{\sigma_{n} \sigma_{n-1} ..\sigma_{1} a_{n+1}}\big)^{k-1} \nonumber\\\nonumber\\
					&=& \big(q-1\big) \big(g_{\sigma_{n} \sigma_{n-1} ..\sigma_{1} a_{n+1}}\big)^{k} + \nonumber\\
					& & q g_{\sigma_{n-1} \sigma_{n-2} ..\sigma_{1}}F\big(t_{a_{n}}\big) g_{\sigma_{n}} \big(g_{\sigma_{n} \sigma_{n-1} ..\sigma_{1} a_{n+1}}\big)^{k-1}  \Psi^{k-1} \big(F((t_{a_{n}})^{-1})\big). \nonumber\\\nonumber
				\end{eqnarray}
							
				So, by induction on $k$, (1) follows. In the very same way we deal with (2), by noticing that: $g_{a_{n+1}}g_{\sigma_{n}} =g^{-1}_{\sigma_{n}} F(t_{a_{n}}) g^{2}_{\sigma_{n}} = (q-1) g_{a_{n+1}} + qg^{-1} _{\sigma_{n}} F(t_{a_{n}}) $.
			\end{demo}
	
 			A main result in \cite{Sadek_2013_2} is to give a general form for ``fully commutative braids'', from which we deduce that any element of the basis of $\widehat{TL}_{n+1} (q)$ (where we have the convention $\sigma_{n+1} = 1$ in $W(\tilde{A_{n}})$ thus $ g_{\sigma_{n} \sigma_{n-1} .. \sigma_{i}} =1$ when $i=n+1$), is either of the form\\
 			$$ c (g_{\sigma_{n} \sigma_{n-1} ..\sigma_{1} a_{n+1}})^{k} g_{\sigma_{n}\sigma_{n-1} .. \sigma_{i}} $$
 			
 			or of the form\\
 			
 			$$g_{\sigma_{i_{0}} .. \sigma_{2}\sigma_{1}a_{n+1}} (g_{\sigma_{n} \sigma_{n-1} ..\sigma_{1} a_{n+1}})^{k} d  g_{\sigma_{n} \sigma_{n-1} .. \sigma_{i}} $$\\
 			
 			where $c$ and $d$ are in $F (\widehat{TL}_{n} (q)) $, $ 1\leq i \leq n+1 $ and $ 0 \leq i_{0} \leq n-1 $ . \\
	
			By Lemma \ref{5_1_3} $ c (g_{\sigma_{n} \sigma_{n-1} ..\sigma_{1} a_{n+1}})^{k} g_{\sigma_{n}\sigma_{n-1} .. \sigma_{i}} $ is of the form:
			\begin{eqnarray}
				\sum^{j=h}_{j=1}c_{j}  (g_{\sigma_{n} \sigma_{n-1} ..\sigma_{1} a_{n+1}})^{j} + M. \nonumber
			\end{eqnarray}
			
			Where $h \leq k $, $c_{j}$ is in  $F (\widehat{TL}_{n} (q)) $ for any $ j $ and $M$ is a Markov element.\\ 
			
			Now we deal with the second form:
			\begin{eqnarray}
				\tau_{n+1} &\big(&g_{\sigma_{i_{0}} .. \sigma_{2}\sigma_{1}a_{n+1}} c (g_{\sigma_{n} \sigma_{n-1} ..\sigma_{1} a_{n+1}})^{k}  g_{\sigma_{n} \sigma_{n-1} .. \sigma_{i}}\big) = \tau_{n+1} \big(g_{\sigma_{n} \sigma_{n-1} .. \sigma_{i}} g_{\sigma_{i_{0}} .. \sigma_{2}\sigma_{1}a_{n+1}} c (g_{\sigma_{n} \sigma_{n-1} ..\sigma_{1} a_{n+1}})^{k}\big). \nonumber
			\end{eqnarray}
	
			For any possible value for $ i_{0} $ or $i$, we see that:
			\begin{eqnarray}
				g_{\sigma_{n} \sigma_{n-1} .. \sigma_{i}} g_{\sigma_{i_{0}} .. \sigma_{2}\sigma_{1}a_{n+1}} c (g_{\sigma_{n} \sigma_{n-1} ..\sigma_{1} a_{n+1}})^{k} = c'g_{\sigma_{n}} (g_{\sigma_{n} \sigma_{n-1} ..\sigma_{1} a_{n+1}})^{s} c'' ,\nonumber
			\end{eqnarray}
			
			where $ c',c''$ are in $F (\widehat{TL}_{n} (q)) $ and $ s \leq k+1 $. By Lemma \ref{5_1_3} we see that this element is of the form:
			\vspace{-0.25cm}  
			\begin{eqnarray}
				\sum^{j=h}_{j=1}f_{j}  (g_{\sigma_{n} \sigma_{n-1} ..\sigma_{1} a_{n+1}})^{j} + M, \nonumber
			\end{eqnarray}
			
			where $h \leq k+1 $, $f_{j}$ is in  $F (\widehat{TL}_{n} (q)) $ for any $ j $ and $M$ is a Markov element .\\ 
	        
			Hence, we see that in order to define $\tau_{n+1}$ uniquely it is enough to have its values on Markov elements and its values on $\Omega (g_{\sigma_{n} \sigma_{n-1} ..\sigma_{1} a_{n+1}})^{k}$, where $1 \leq k$ (since if $k$ is equal to 0 then we are again in the case of a Markov element) and $\Omega$ is in $F\big(\widehat{TL}_{n} (q)\big)$.
			\vspace{0.25cm}
			\begin{lemme}\label{5_1_4}
				Let $2 \leq n $ then  $\tau_{n+1}$ is uniquely defined by its values on Markov elements, in addition to its values on $ (g_{\sigma_{n} \sigma_{n-1} ..\sigma_{1} a_{n+1}} )^{k} $,  with $0 \leq k $ .\\
			\end{lemme}

			\begin{demo} 
		 		In order to determine $\tau_{n+1} \big( h(g_{\sigma_{n} \sigma_{n-1} ..\sigma_{1} a_{n+1}})^{k}\big)$, with a positive $k$ and an arbitrary $h$ in  $F\big(\widehat{TL}_{n} (q)\big)$, it is enough to treat $\tau_{n+1} \big(F(t_{x})(g_{\sigma_{n} \sigma_{n-1} ..\sigma_{1} a_{n+1}})^{k}\big)$, with $x$ in $W^{c}(\tilde{A_{n-1}})$, but the fact that $\tau_{n+1}$ is a trace, in addition to the fact that $g_{\sigma_{n} \sigma_{n-1} ..\sigma_{1} a_{n+1}}$ acts as a Dynkin automorphism on $F \big(\widehat{TL}_{n} (q)\big)$, authorizes us to suppose that $x$ has a reduced expression which ends with $\sigma_{n-1}$. \\
		 		
		 		Now we show by induction on $l(x)$, that $\tau_{n+1} \big( F(t_{x}) (g_{\sigma_{n} \sigma_{n-1} ..\sigma_{1} a_{n+1}})^{k}\big)$ is a sum of values of $\tau_{n+1}$ over $ (g_{\sigma_{n} \sigma_{n-1} ..\sigma_{1} a_{n+1}})^{k}$, elements of the form $ h (g_{\sigma_{n} \sigma_{n-1} ..\sigma_{1} a_{n+1}})^{i}$ with $i < k $ and Markov elements, (of course with coefficients in the ground ring which might be zeros).\\ 
	
				For $l(x) = 0$ the property is true. Take $l(x) > 0$, and let $ x = z \sigma_{n-1}$ be a reduced expression, hence:
				\begin{eqnarray}
					\tau_{n+1} \big( F(t_{x}) (g_{\sigma_{n} \sigma_{n-1} ..\sigma_{1} a_{n+1}})^{k}\big) &=& \tau_{n+1}\big( F(t_{z}) F(t_{\sigma_{n-1}}) g_{\sigma_{n} \sigma_{n-1} ..\sigma_{1} a_{n+1}}(g_{\sigma_{n} \sigma_{n-1} ..\sigma_{1} a_{n+1}})^{k-1}\big) \nonumber\\\nonumber\\
					&=& \tau_{n+1} \big(F(t_{z}) \underbrace{g_{\sigma_{n-1}} g_{\sigma_{n}} g_{\sigma_{n-1}}}_{=-V(g_{\sigma_{n-1}}, g_{\sigma_{n}})} g_{ \sigma_{n-2} ..\sigma_{1} a_{n+1}}(g_{\sigma_{n} \sigma_{n-1} ..\sigma_{1} a_{n+1}})^{k-1}\big). \nonumber
				\end{eqnarray}
	
				Recalling that $V(g_{\sigma_{n-1}}, g_{\sigma_{n}})=0$, this is equal to the following sum:
				\begin{eqnarray}
					& & - \tau_{n+1}\big( F(t_{z}) (g_{\sigma_{n} \sigma_{n-1} ..\sigma_{1} a_{n+1}})^{k}\big) \nonumber\\\nonumber\\
					& & - \tau_{n+1}\big( F(t_{z}) g_{\sigma_{n-1}}   g_{ \sigma_{n-2} ..\sigma_{1}}  g_{a_{n+1}}(g_{\sigma_{n} \sigma_{n-1} ..\sigma_{1} a_{n+1}})^{k-1}\big)\nonumber\\\nonumber\\
					& & - \tau_{n+1}\big( F(t_{z}) g_{ \sigma_{n-2} ..\sigma_{1} a_{n+1}}(g_{\sigma_{n} \sigma_{n-1} ..\sigma_{1} a_{n+1}})^{k-1}\big)\nonumber\\\nonumber\\
					& & - \tau_{n+1}\big( F(t_{z}) g_{\sigma_{n-1}}   g_{ \sigma_{n-2} ..\sigma_{1}}\underbrace{ g_{\sigma_{n}} g_{a_{n+1}}}_{}(g_{\sigma_{n} \sigma_{n-1} ..\sigma_{1} a_{n+1}})^{k-1}\big) \nonumber\\
					& & - \tau_{n+1}\big( F(t_{z}) g_{ \sigma_{n-2} ..\sigma_{1}} \underbrace{g_{\sigma_{n}} g_{a_{n+1}}}_{}(g_{\sigma_{n} \sigma_{n-1} ..	\sigma_{1} a_{n+1}})^{k-1}\big).\nonumber
				\end{eqnarray}
				
				Now we apply the induction hypothesis to the first term. The second and the third terms are equal to:
				\begin{eqnarray}
					\tau_{n+1} &\bigg(& F\big(t_{z}\big) g_{\sigma_{n-1}} g_{ \sigma_{n-2} ..\sigma_{1}} F\big( t_{a_{n}}\big) g_{\sigma_{n}} F\big( (t_{a_{n}})^{-1}\big)\big(g_{\sigma_{n} \sigma_{n-1} ..\sigma_{1} a_{n+1}}\big)^{k-1}\bigg)\nonumber\\\nonumber\\
					&+& \tau_{n+1} \bigg( F\big(t_{z}\big) g_{\sigma_{n-2} ..\sigma_{1}} F\big( t_{a_{n}}\big) g_{\sigma_{n}}F\big( (t_{a_{n}})^{-1}\big)\big(g_{\sigma_{n} \sigma_{n-1} ..\sigma_{1} a_{n+1}}\big)^{k-1}\bigg), \nonumber					
				\end{eqnarray}
				
				which is equal to:
				\begin{eqnarray}
					\tau_{n+1} &\bigg(& \Psi^{1-k} \big( F\big( (t_{a_{n}})^{-1}\big) \big)F\big(t_{z}\big) g_{\sigma_{n-1}}   g_{ \sigma_{n-2} ..\sigma_{1}} F\big( t_{a_{n}}\big) \big(g_{\sigma_{n}} (g_{\sigma_{n} \sigma_{n-1} ..\sigma_{1} a_{n+1}})^{k-1}\big)\bigg)\nonumber\\\nonumber\\	
					&+& \tau_{n+1}\bigg(\Psi^{1-k} \big( F\big( (t_{a_{n}})^{-1}\big) \big) F\big(t_{z}\big) g_{ \sigma_{n-2} ..\sigma_{1}} F( t_{a_{n}})\big( g_{\sigma_{n}}(g_{\sigma_{n} \sigma_{n-1} ..\sigma_{1} a_{n+1}})^{k-1}\big)\bigg).\nonumber				
				\end{eqnarray}
				
				The fourth and the fifth terms are equal to:
				\begin{eqnarray}
					\tau_{n+1} &\bigg(& F\big(t_{z}\big) g_{\sigma_{n-1}} g_{ \sigma_{n-2} ..\sigma_{1}} F\big( t_{a_{n}}\big) \big(g_{\sigma_{n}}(g_{\sigma_{n} \sigma_{n-1} ..\sigma_{1} a_{n+1}})^{k-1}\big)\bigg) \nonumber\\\nonumber\\
					&+& \tau_{n+1}\bigg( F\big(t_{z}\big) g_{\sigma_{n-2} ..\sigma_{1}} F\big(t_{a_{n}}\big) \big(g_{\sigma_{n}}(g_{\sigma_{n} \sigma_{n-1} ..\sigma_{1} a_{n+1}})^{k-1}\big)\bigg).\nonumber
				\end{eqnarray}
	
				Thus, Lemma \ref{5_1_3} tells us that the property is true for those four terms. This step is to be applied repeatedly, to the powers of $g_{\sigma_{n} \sigma_{n-1} ..\sigma_{1} a_{n+1}} $ down to an element of the form $ \tau_{n+1}\big(h (g_{\sigma_{n} \sigma_{n-1} ..\sigma_{1} a_{n+1}})^{1}\big)$, arriving to the sum of:	
				$$\tau_{n+1}( g_{\sigma_{n} \sigma_{n-1} ..\sigma_{1} a_{n+1}})$$ and $$\tau_{n+1} (h' g_{\sigma_{n-1} ..\sigma_{1} a_{n+1}}) ,$$ 				
				which is the sum of  values of $ \tau_{n+1}$ on Markov elements, since $h,h'\in F \big(\widehat{TL}_{n} (q)\big)$.
	\end{demo}
			We end this part by the following Lemma:			
			
			\begin{lemme} \label{5_1_5}
				Let $ 1\leq k $. Then $ (g_{\sigma_{n} \sigma_{n-1} ..\sigma_{1} a_{n+1}})^{k}$ is a sum of two kinds of elements:\\
				\begin{itemize} 
					\item[$(1)$] $g_{\sigma_{n}} \big(g_{\sigma_{n-1} \sigma_{n-2} .. \sigma_{1}}  F(t_{a_{n}}) \big)^{j} g_{\sigma_{n}} h $, with $j\leq k $.
					\item[$(2)$] $\big(g_{\sigma_{n-1} \sigma_{n-2} .. \sigma_{1}}  F(t_{a_{n}})\big)^{i} g_{\sigma_{n}} f $, with $i < k $,\\
				\end{itemize}	
					
				with $h,f$ in $F \big(\widehat{TL}_{n} (q)\big)$ and $2\leq n$.\\
    
				Moreover, in the first type we have one, and only one element, with $j=k$, in which we have: 
				\vspace{-0.25cm}
				\begin{eqnarray}
					h = \prod^{i=k-1}_{i=0} \Psi^{-i}\big(F(t^{-1}_{a_{n}})\big). \nonumber
				\end{eqnarray}
			\end{lemme}
      
			\begin{demo}
				Suppose that $k=1$. Then,
				\begin{eqnarray}
					g_{\sigma_{n} \sigma_{n-1} ..\sigma_{1} a_{n+1}} = g_{\sigma_{n}} \big(g_{\sigma_{n-1} \sigma_{n-2} .. \sigma_{1}} F(t_{a_{n}})\big) g_{\sigma_{n}} F\big(t_{a_{n}}\big)^{-1}, \nonumber
				\end{eqnarray}
				
				the property is true.\\ 

				Suppose the property is true for $k-1$, then, with $ 2\leq k$, we have:
				\begin{eqnarray}
					(g_{\sigma_{n} \sigma_{n-1} ..\sigma_{1} a_{n+1}})^{k} = (g_{\sigma_{n} \sigma_{n-1} ..\sigma_{1} a_{n+1}})^{k-1} g_{\sigma_{n} \sigma_{n-1} ..\sigma_{1} a_{n+1}}. \nonumber
				\end{eqnarray}
				
				We apply the property to $(g_{\sigma_{n} \sigma_{n-1} ..\sigma_{1} a_{n+1}})^{k-1}$, which gives two cases:\\
				
				\begin{itemize} 
					\item[$(1)$] $g_{\sigma_{n}} \big(g_{\sigma_{n-1} \sigma_{n-2} .. \sigma_{1}}  F(t_{a_{n}}) \big)^{j'} g_{\sigma_{n}} h  g_{\sigma_{n} \sigma_{n-1} ..\sigma_{1} a_{n+1}}$, with $ j' \leq k-1 $ which is: \textcolor{white}{...................} $g_{\sigma_{n}} \big(g_{\sigma_{n-1} \sigma_{n-2} .. \sigma_{1}}  F(t_{a_{n}}) \big)^{j'} g_{\sigma_{n}} g_{\sigma_{n} \sigma_{n-1} ..\sigma_{1} a_{n+1}} \Psi^{-1}\big(h\big)$, which is equal to:\\
					\begin{eqnarray}
						qg_{\sigma_{n}} &\big(& g_{\sigma_{n-1} \sigma_{n-2} .. \sigma_{1}}  F(t_{a_{n}}) \big)^{j'+1} g_{\sigma_{n}} F\big( (t_{a_{n}})^{-1}\big) \Psi^{-1}\big(h\big) \nonumber\\\nonumber\\
						&+& (q-1) g_{\sigma_{n}}  g_{\sigma_{n} \sigma_{n-1} ..\sigma_{1} a_{n+1}}\Psi^{-1} \bigg(\big(g_{\sigma_{n-1} \sigma_{n-2} .. \sigma_{1}}  F(t_{a_{n}}) \big)^{j'} \bigg)   \Psi^{-1}\left(h\right). \nonumber
					\end{eqnarray}

					Since, $ j'+1 \leq k $, the first term is clear to be of the first type, while the second term is equal to:
					\begin{eqnarray}
						& & \big(q-1\big)q g_{\sigma_{n-1} ..\sigma_{1}} F\big(t_{a_{n}}\big) g_{\sigma_{n}} F\big( (t_{a_{n}})^{-1}\big)\Psi^{-1} \big(\big(g_{\sigma_{n-1} \sigma_{n-2} .. \sigma_{1}}  F(t_{a_{n}}) \big)^{j'} \big)   \Psi^{-1}\big(h\big) + \nonumber\\\nonumber\\
						& & \big(q-1\big)^{2} g_{\sigma_{n} \sigma_{n-1} ..\sigma_{1} a_{n+1}}\Psi^{-1} \big(\big(g_{\sigma_{n-1} \sigma_{n-2} .. \sigma_{1}}  F(t_{a_{n}}) \big)^{j'} \big) \Psi^{-1}\big(h\big) .\nonumber
					\end{eqnarray}

					Here, the first term is of the second type (with $i=1 < k $), and the second term is of the first type (with $j=1$).\\
					
					\item[$(2)$] $\big(g_{\sigma_{n-1} \sigma_{n-2} .. \sigma_{1}}  F(t_{a_{n}})\big)^{i'} g_{\sigma_{n}} f g_{\sigma_{n} \sigma_{n-1} ..\sigma_{1} a_{n+1}}$, with $ i' < k-1 $, which is:
					\begin{eqnarray}
						\big(&g&_{\sigma_{n-1} \sigma_{n-2} .. \sigma_{1}}  F(t_{a_{n}}) \big)^{i'} g_{\sigma_{n}}  g_{\sigma_{n} \sigma_{n-1} ..\sigma_{1} a_{n+1}} \Psi^{-1} \big(f\big) = \nonumber\\\nonumber\\
						& & q \big(g_{\sigma_{n-1} \sigma_{n-2} .. \sigma_{1}}  F(t_{a_{n}})\big)^{i'+1} g_{\sigma_{n}} F\big((t_{a_{n}})^{-1}\big) \Psi^{-1} \big(f\big) + \nonumber\\\nonumber\\
						& & \big(q-1\big) g_{\sigma_{n}} \big(g_{\sigma_{n-1} ..\sigma_{1}}  F(t_{a_{n}})\big) g_{\sigma_{n}}  F\big((t_{a_{n}})^{-1}\big) \Psi^{-1}\big(\big(g_{\sigma_{n-1} \sigma_{n-2} .. \sigma_{1}}  F(t_{a_{n}}) \big)^{i'} \big)  \Psi^{-1} \big(f\big).  \nonumber\\\nonumber
					\end{eqnarray}
					
					Since $i'+1 <k $, the first term is of the second type, while the second term is of the first type with $j=1$. The Lemma is proven. \\

					(By induction over $k$ again, the last formula is easy). 
				\end{itemize}
				
			\end{demo}

			\textbf{Part 2}\\
			
			In this part we treat Theorem   \ref{5_1_1} when $ n\geq 3 $. As said at the beginning of Part 1, for $ n \geq 3 $, and by sending the ``fully commutative braids'' onto $\widehat{TL}_{n} (q)$, we get that any element of the basis of $ \widehat{TL}_{n+1}(q) $ is a linear combination of two kinds of elements, namely:
			\begin{eqnarray}
				I &=& F_{n}(t_{u}) (g_{\sigma_{n} \sigma_{n-1} ..\sigma_{1} a_{n+1}})^{k} g_{\sigma_{n}\sigma_{n-1} .. \sigma_{s}}, \nonumber\\\nonumber\\
				II &=& g_{\sigma_{i_{0}} .. \sigma_{2}\sigma_{1}a_{n+1}} F_{n}(t_{u}) (g_{\sigma_{n} \sigma_{n-1} ..\sigma_{1} a_{n+1}})^{k} g_{\sigma_{n}\sigma_{n-1} .. \sigma_{s}}, \nonumber
			\end{eqnarray}
			
			here, $u$ is in $W^{c}(\tilde{A_{n-1}})$, where $1 \leq s \leq n+1$ with $ 0 \leq i_{0} \leq n-1 $ and $ 0\leq k$.\\
			
			Using Lemma \ref{5_1_3} we see that:
			\begin{eqnarray}
				\big( g_{\sigma_{n} \sigma_{n-1} ..\sigma_{1} a_{n+1}} \big)^{k}g_{\sigma_{n}} &=& (q-1) \big(g_{\sigma_{n} \sigma_{n-1} ..\sigma_{1} a_{n+1}} \big)^{k} \nonumber\\\nonumber\\
				& & + \sum^{i=k-1}_{i=1} h_{i} \big( g_{\sigma_{n} \sigma_{n-1} ..\sigma_{1} a_{n+1}} \big)^{i} \nonumber\\\nonumber\\
				& & + A \prod^{j=k-1}_{j=0} \Psi^{-j} \big(\big(g_{\sigma_{{n-1}}} \big)^{-1} \big) g_{\sigma_{n}} \big(g_{\sigma_{n-1} \sigma_{n-2} ..\sigma_{1}} F(t_{a_{n}}) \big)^{k}, \nonumber
			\end{eqnarray}
			
			but, $I = F_{n} (t_{u}) \underbrace{(g_{\sigma_{n} \sigma_{n-1} ..\sigma_{1} a_{n+1}})^{k} g_{\sigma_{n}}}_{}g_{\sigma_{n-1} .. \sigma_{s}} $, that is:
			\begin{eqnarray}
				I &=& (q-1) F_{n}\big(t_{u}\big) \big(g_{\sigma_{n} \sigma_{n-1} ..\sigma_{1} a_{n+1}} \big)^{k} g_{\sigma_{n-1} .. \sigma_{s}} \nonumber\\\nonumber\\
				& & + \sum^{i=k-1}_{i=1} h_{i} F_{n}\big(t_{u}\big) \big(g_{\sigma_{n} \sigma_{n-1} ..\sigma_{1} a_{n+1}}\big)^{i} g_{\sigma_{n-1} .. \sigma_{s}} \nonumber\\\nonumber\\
				& & + A \prod^{j=k-1}_{j=0} F_{n}\big(t_{u}\big) \Psi^{-j} \big( \big( g_{\sigma_{{n-1}}} \big)^{-1}\big) g_{\sigma_{n}} \big(g_{\sigma_{n-1} \sigma_{n-2} ..\sigma_{1}}F(t_{a_{n}}) \big)^{k} g_{\sigma_{n-1} .. \sigma_{s}}. \nonumber
			\end{eqnarray}

			Using the action of $g_{\sigma_{n} \sigma_{n-1} ..\sigma_{1} a_{n+1}} $ on $ F_{n}(\widehat{TL}_{n}(q))$, we see that:
			\begin{eqnarray}
				I = \sum^{i=k}_{i=1} F_{n} \big(t_{b_{i}} \big) \big(g_{\sigma_{n} \sigma_{n-1} ..\sigma_{1} a_{n+1}} \big)^{i} + \sum_{j} F_{n} \big(t_{b_{j}} \big) g_{\sigma_{n}} F_{n} \big(t_{d_{j}} \big), \nonumber
			\end{eqnarray}
			
			where $ b_{j} $, $c_{j}$ and $ d_{i} $ are in $W^{c}(\tilde{A_{n-1}})$, for every $i$ and $j$.\\

			Now, we see, as well, that:  
			\begin{eqnarray}
				II &=& g_{\sigma_{i_{0}} .. \sigma_{2} \sigma_{1}} g_{a_{n+1}} F_{n} \big(t_{u}\big) \big( g_{\sigma_{n} \sigma_{n-1} ..\sigma_{1} a_{n+1}} \big)^{k} g_{\sigma_{n}\sigma_{n-1} .. \sigma_{s}} \nonumber\\\nonumber\\
				&=& g_{\sigma_{i_{0}} .. \sigma_{2} \sigma_{1}} F_{n} \big(t_{a_{n}}\big)\underbrace{g_{\sigma_{n}} \big(g_{\sigma_{n} \sigma_{n-1} ..\sigma_{1} a_{n+1}} \big)^{k}}_{}\Psi^{k} \big( F_{n} \big(t^{-1}_{a_{n}} \big) F_{n}\big(t_{u}\big) \big) g_{\sigma_{n}\sigma_{n-1} .. \sigma_{s}}, \nonumber
			\end{eqnarray}
			
			since $g_{a_{n+1}} = F_{n}(t_{a_{n}}) g_{\sigma_{n}} F_{n}(t^{-1}_{a_{n}})$.\\ 

			By Lemma \ref{5_1_3}, we see that $ II $ is equal to:
			\begin{eqnarray}
				& & (q-1) g_{\sigma_{i_{0}} .. \sigma_{2}\sigma_{1}} F_{n} \big(t_{a_{n}} \big) \big(g_{\sigma_{n} \sigma_{n-1} ..\sigma_{1} a_{n+1}} \big)^{k} \Psi^{k} \big(F_{n} \big( t^{-1}_{a_{n}} \big) F_{n}(t_{u}) \big) g_{\sigma_{n}\sigma_{n-1} .. \sigma_{s}} \nonumber\\\nonumber\\
				& & + \sum^{i=k-1}_{i=1} g_{\sigma_{i_{0}} .. \sigma_{2}\sigma_{1}} F_{n} \big(t_{a_{n}}\big) f_{i} \big(g_{\sigma_{n} \sigma_{n-1} ..\sigma_{1} a_{n+1}} \big)^{i} \Psi^{k} \big(F_{n}(t^{-1}_{a_{n}})F_{n}(t_{u}) \big) g_{\sigma_{n}\sigma_{n-1} .. \sigma_{s}} \nonumber\\\nonumber\\
				& & + Ag_{\sigma_{i_{0}} .. \sigma_{2}\sigma_{1}} F_{n} \big(t_{a_{n}}\big) \big(g_{\sigma_{n-1} \sigma_{n-2} ..\sigma_{1}} F(t_{a_{n}}) \big)^{k} g_{\sigma_{n}} \prod^{j=k-1}_{j=0} \Psi^{j} \big( F\big((t_{a_{n}})^{-1}\big) \big) \Psi^{k} \big(F_{n}( t^{-1}_{a_{n}}) F_{n}(t_{u}) \big) g_{\sigma_{n} \sigma_{n-1} .. \sigma_{s}}, \nonumber 
			\end{eqnarray}
			
			which is equal to:
			\begin{eqnarray}
				& &  \sum^{i=k}_{i=1} F_{n} \big(t_{x'_{i}}\big) \big(\underbrace{g_{\sigma_{n} \sigma_{n-1} ..\sigma_{1} a_{n+1}} \big)^{i} g_{\sigma_{n}}}_{}g_{\sigma_{n-1} .. \sigma_{s}}  \nonumber\\\nonumber\\
				& & + Ag_{\sigma_{i_{0}} .. \sigma_{2} \sigma_{1}} F_{n} \big(t_{a_{n}}\big) \big(g_{\sigma_{n-1} \sigma_{n-2} ..\sigma_{1}} F(t_{a_{n}}) \big)^{k} g_{\sigma_{n}} \prod^{j=k-1}_{j=0} \Psi^{j} \big( F\big((t_{a_{n}})^{-1}\big) \big) \Psi^{k} \big(F_{n}(t^{-1}_{a_{n}})F_{n}(t_{u}) \big) g_{\sigma_{n}	\sigma_{n-1} .. \sigma_{s}}, \nonumber 
			\end{eqnarray}
	
			where $ x'_{i} $ is in $W^{c}(\tilde{A_{n-1}})$ for all $i$. \\ 
	
			Now we repeat the same step as for $I$, to get the next corollary.  \\ 
	
			\begin{corollaire} \label{5_1_9}
				Let $ 3 \leq n$. Let $w$ be in $W^{c}(\tilde{A_{n}})$ . \\ 
	
				Then there exist $ 0 \leq k $ and  $ 1 \leq s \leq n+1 $, there exist $x_{i}$ , $y_{j}$ and $ z_{j}$ in $W^{c}(\tilde{A_{n-1}})$ $\big[$With the convention $W^{c}(\tilde{A_{2}}) = W(\tilde{A_{2}}) \big]$ such that:
				\begin{eqnarray}			
					g_{w} = \sum^{i=k}_{i=1} F_{n}\big(t_{x_{i}}\big) \big(g_{\sigma_{n} \sigma_{n-1} ..\sigma_{1} a_{n+1}} \big)^{i} + \sum_{j} F_{n} \big(t_{y_{j}} \big) g_{\sigma_{n}} F_{n}\big(t_{z_{j}}\big) g_{\sigma_{n}\sigma_{n-1} .. \sigma_{s}}. \nonumber 
				\end{eqnarray}
			\end{corollaire} 
	
			Now we suppose that $ 3 \leq n $. Consider the following sequence:\\
			\begin{eqnarray}
				\widehat{TL}_{n-1}(q) \longrightarrow \widehat{TL}_{n}(q) \longrightarrow  \widehat{TL}_{n+1}(q). \nonumber 
			\end{eqnarray}
			
			We keep using $t_{\sigma_{i}}$ (resp. $g_{\sigma_{i}}$) as generators of  $\widehat{TL}_{n}(q) $ (resp. $\widehat{TL}_{n+1}(q) $). We use $e_{\sigma_{i}}$ for $\widehat{TL}_{n-1}(q)$. With a simple computation, we see that $g_{\sigma_{n}}$ commutes with $ F_{n}F_{n-1}(e_{\sigma_{i}}) $, for $1 \leq i \leq n-2$, and with $ F_{n}F_{n-1}(e_{a_{n-1}}) $, hence it commutes with every element in $ F_{n}F_{n-1}(\widehat{TL}_{n-1}(q))$.  \\

			Lemma \ref{5_1_4} and Lemma \ref{5_1_5} confirm that $\tau_{n+1} $ is uniquely defined over $\widehat{TL}_{n+1}(q) $ by its values on  $g_{\sigma_{n}} (g_{\sigma_{n-1} \sigma_{n-2} .. \sigma_{1}}  F(t_{a_{n}}) )^{k} g_{\sigma_{n}} h $, for a positive $ k $ and an arbitrary $h$ in $F (\widehat{TL}_{n} (q)) $ beside its values on Markov elements. In other terms:  $\tau_{n+1} $ is uniquely defined over $\widehat{TL}_{n+1}(q) $ by its values over $g_{\sigma_{n}} \big( g_{\sigma_{n-1} \sigma_{n-2} .. \sigma_{1}}  F(t_{a_{n}}) \big)^{k} g_{\sigma_{n}} F_{n} \big(t_{v}\big)$, with a positive $ k $ and an arbitrary $v$ in $ W^{c}(\tilde{A_{n-1}}) $, besides the values on Markov elements.
			\begin{eqnarray}
				\text{Set }I := \tau_{n+1} \bigg( g_{\sigma_{n}} \big( g_{\sigma_{n-1} \sigma_{n-2} .. \sigma_{1}}  F(t_{a_{n}}) \big)^{k} g_{\sigma_{n}} F_{n} \big(t_{v} \big)\bigg),  \nonumber
			\end{eqnarray}			
			
			by Corollary \ref{5_1_9} we see that:
			\begin{eqnarray}
				t_{v} &=& \sum^{i=h}_{i=1} \underbrace{ F_{n-1} \big(e_{x_{i}} \big) \big(t_{\sigma_{n-1} \sigma_{n-2} ..\sigma_{1} a_{n}} \big)^{i}}_{C} \nonumber\\\nonumber\\
				&+& \sum_{j} \underbrace{ F_{n-1} \big(e_{y_{j}} \big) t_{\sigma_{n-1}} F_{n-1} \big(e_{z_{j}} \big) t_{\sigma_{n-1}\sigma_{n-2} .. \sigma_{s}}}_{B} \nonumber\\\nonumber\\
				&+& \sum_{j} \underbrace{ F_{n-1} \big(e_{y'_{j}} \big) t_{\sigma_{n-1}} F_{n-1} \big(e_{z'_{j}} \big)}_{A} ,\nonumber\\\nonumber					
			\end{eqnarray}
			
			where $ 0 \leq h $ and  $ 1 \leq s \leq n-1 $. With $ x_{i}$ , $y_{j}$, $ z_{j}$, $y'_{j}$ and $ z'_{j}$ are in $W^{c}(\tilde{A_{n-2}})$.\\

			We have added the third term $C$ to the two terms of Corollary \ref{5_1_9}, because we had to take into account here, the case of  $s= n+1 $, i.e.,  $ g_{\sigma_{n+1}} = 1$ for $ W^{c}(\tilde{A_{n-1}}) $.\\
 
			For terms of \textbf{Type (A)}, we see that:
			\begin{eqnarray}
 				I_{1} &:=& \tau_{n+1} \bigg(g_{\sigma_{n}} \big(g_{\sigma_{n-1} \sigma_{n-2} .. \sigma_{1}}  F_{n}(t_{a_{n}}) \big)^{k} g_{\sigma_{n}} F_{n} \big(F_{  n -1}(e_{y'_{j}}) t_{\sigma_{n-1}} F_{n-1}(e_{z'_{j}}) \big) \bigg) \nonumber\\\nonumber\\
 				&=& \tau_{n+1} \bigg( \big(g_{\sigma_{n-1} \sigma_{n-2} .. \sigma_{1}} F_{n}(t_{a_{n}}) \big)^{k} F_{n} \big(F_{n-1}(e_{y'_{j}}) \big) g_{\sigma_{n}} F_{n} \big(t_{\sigma_{n-1}} \big) g_{\sigma_{n}} F_{n} \big( F_{n-1}(e_{z'_{j}}) \big) \bigg) \nonumber
 			\end{eqnarray}	
 			
 			\begin{eqnarray}	
 				&=& \tau_{n+1} \bigg( \big( g_{\sigma_{n-1} \sigma_{n-2} .. \sigma_{1}} F(t_{a_{n}}) \big)^{k} F_{n} \big(F_{n-1}(e_{y'_{j}}) \big) \underbrace{ g_{\sigma_{n}} g_{\sigma_{n-1}} g_{\sigma_{n}}}_{r}   F_{n} \big( F_{n-1}(e_{z'_{j}})\big) \bigg) \nonumber\\\nonumber\\
 				&=& \tau_{n+1} \bigg( \big( g_{\sigma_{n-1} \sigma_{n-2} .. \sigma_{1}}  F(t_{a_{n}}) \big)^{k} F_{n}\big(F_{n-1}(e_{y'_{j}}) \big) \underbrace{ g_{\sigma_{n-1}}  g_{\sigma_{n}} g_{\sigma_{n-1}}}_{r}   F_{n} \big( F_{n-1}(e_{z'_{j}}) \big) \bigg), \nonumber			
			\end{eqnarray}

			which is clearly, the sum of values of $\tau_{n+1} $ on Markov elements.\\ 

			For terms of \textbf{Type (B)}, we see that:
			\begin{eqnarray}
				I_{2} &:=& \tau_{n+1} \bigg[ g_{\sigma_{n}} \big(g_{\sigma_{n-1} \sigma_{n-2} .. \sigma_{1}}  F(t_{a_{n}}) \big)^{k} g_{\sigma_{n}} F_{n} \big( F_{n-1}(e_{y_{j}}) t_{\sigma_{n-1}} F_{n-1}(e_{z_{j}}) t_{\sigma_{n-1}\sigma_{n-2} .. \sigma_{s}} \big) \bigg] \nonumber\\\nonumber\\
				&=& \tau_{n+1} \bigg[ g_{\sigma_{n}} \big( g_{\sigma_{n-1} \sigma_{n-2} .. \sigma_{1}}  F(t_{a_{n}}) \big)^{k} g_{\sigma_{n}} F_{n} \big( F_{n-1}(e_{y_{j}}) t_{\sigma_{n-1}} F_{n-1}(e_{z_{j}}) t_{\sigma_{n-1}}t_{\sigma_{n-2} .. \sigma_{s}} \big) \bigg] \nonumber\\\nonumber\\
				&=& \tau_{n+1} \bigg[ g_{\sigma_{n}} F_{n} F_{n-1} \big(e_{\sigma_{n-2} .. \sigma_{s}} \big) \big(g_{\sigma_{n-1} \sigma_{n-2} .. \sigma_{1}}  F(t_{a_{n}}) \big)^{k} g_{\sigma_{n}} F_{n} \big( F_{n-1}(e_{y_{j}}) t_{\sigma_{n-1}} F_{n-1}(e_{z_{j}}) t_{\sigma_{n-1}} \big) \bigg]. \nonumber			
			\end{eqnarray}

			Now, we set $ ^{m}_{r}F:= F_{m}F_{m-1} .. F_{r} $.\\ 

			We call $\delta$ the image of $  F_{n-1}(e_{\sigma_{n-2} .. \sigma_{s}})$ under the action of  $\big(g_{\sigma_{n-1} \sigma_{n-2} .. \sigma_{1}}  F(t_{a_{n}}) \big)^{k}$, thus:
			\begin{eqnarray}
				I_{2} &=& \tau_{n+1} \bigg[ g_{\sigma_{n}} \big( g_{\sigma_{n-1} \sigma_{n-2} .. \sigma_{1}}  F(t_{a_{n}}) \big) ^{k} g_{\sigma_{n}} F_{n} \big(\delta\big) \big(^{n}_{n-1}F(e_{y_{j}}) \big) g_{\sigma_{n-1}}  \big(^{n}_{n-1}F(e_{z_{j}}) \big) g_{\sigma_{n-1}} \bigg] \nonumber\\\nonumber\\
				&=& \tau_{n+1} \bigg[ g_{\sigma_{n}} \big( F_{n} (t_{\sigma_{n-1} \sigma_{n-2} .. \sigma_{1}a_{n}}) \big) ^{k} g_{\sigma_{n}} F_{n} \big(\delta \big) \big(^{n}_{n-1}F(e_{y_{j}}) \big) g_{\sigma_{n-1}}  \big(^{n}_{n-1}F(e_{z_{j}}) \big) g_{\sigma_{n-1}} \bigg].\nonumber			
			\end{eqnarray}
			
			Now consider $ (t_{\sigma_{n-1} \sigma_{n-2} .. \sigma_{1}a_{n} } )^{k}$. We apply Lemma \ref{5_1_5} to this element in $\widehat{TL}_{n}(q) $, hence, it is the sum of two kind of elements: (1) Markov elements (2) elements of the form $t_{\sigma_{n-1}} (e_{\sigma_{n-2} \sigma_{n-3} .. \sigma_{1}a_{n-1}} )^{j}t_{\sigma_{n-1}} \delta$, where $j \leq k $, and $\delta $ in $F_{n-1}\big(\widehat{TL}_{n-1}(q)\big)$. In the case (1) we are done. If we are in case (2), then we apply the Lemma \ref{5_1_5} on $(e_{\sigma_{n-2} \sigma_{n-3} .. \sigma_{1}a_{n-1}} )^{j} $. We keep going in the same manner, by applying Lemma \ref{5_1_5} repeatedly (in fact $ n-2 $ times), we arrive to:
			\begin{eqnarray}
				t_{\sigma_{n-1}} t_{\sigma_{n-2}} .. t_{\sigma_{2}} &\big(&F_{n-1}F_{n-2} .. F_{2}(^{2} g_{\sigma_{1}a_{2}})^{j} \big) t_{\sigma_{2}} .. t_{\sigma_{n-2}} t_{\sigma_{n-1}} \lambda \nonumber\\\nonumber\\
				&=& t_{\sigma_{n-1}} t_{\sigma_{n-2}} .. t_{\sigma_{2}} \big(^{n-1}_{2} F(^{2} g_{\sigma_{1}a_{2}})^{j} \big) t_{\sigma_{2}} ..  t_{\sigma_{n-2}} t_{\sigma_{n-1}} \lambda, \nonumber			
			\end{eqnarray}

			where $\lambda $ is in $^{n}_{n-1}F(\widehat{TL}_{n-1}(q))$. Leaving aside Markov elements, we get: 
			\begin{eqnarray}
				I_{2} &=& \tau_{n+1} \bigg[ g_{\sigma_{n}} F_{n} \bigg( t_{\sigma_{n-1}} t_{\sigma_{n-2}} ..  t_{\sigma_{2}} \big(^{n-1}_{2}F(^{2}g_{\sigma_{1}a_{2}})^{j} \big) t_{\sigma_{2}} .. t_{\sigma_{n-2}} t_{\sigma_{n-1}} \lambda \bigg) g_{\sigma_{n}} F_{n} \big(\delta \big) \nonumber\\\nonumber\\
				& & \big(^{n}_{n-1}F(e_{y_{j}}) \big) g_{\sigma_{n-1}} \big(^{n}_{n-1}F(e_{z_{j}}) \big) g_{\sigma_{n-1}} \big] \nonumber\\\nonumber\\
				&=& \tau_{n+1} \bigg[ g_{\sigma_{n}} g_{\sigma_{n-1}} g_{\sigma_{n-2}} .. g_{\sigma_{2}} \big(^{n-1}_{2} F(^{2} g_{\sigma_{1}a_{2}})^{j} \big) g_{\sigma_{2}} .. g_{\sigma_{n-2}} g_{\sigma_{n-1}} g_{\sigma_{n}} F_{n} (\lambda \delta) \nonumber\\\nonumber\\
				& & \big(^{n}_{n-1}F(e_{y_{j}}) \big) g_{\sigma_{n-1}} \big(^{n}_{n-1}F(e_{z_{j}}) \big) g_{\sigma_{n-1}} \bigg].\nonumber
			\end{eqnarray}
			
			We set $ M':= F_{n} \big(\lambda \delta \big) \big(^{n}_{n-1}F(e_{y_{j}}) \big) g_{\sigma_{n-1}}  \big(^{n}_{n-1}F(e_{z_{j}}) \big)$, which is a Markov element in $\widehat{TL}_{n-1}(q) $. Hence, we have:
			\begin{eqnarray}
				& & \tau_{n+1} \bigg[ g_{\sigma_{n}} g_{\sigma_{n-1}} g_{\sigma_{n-2}} .. g_{\sigma_{2}} \big(^{n-1}_{2}F(^{2}g_{\sigma_{1}a_{2}})^{j} \big) g_{\sigma_{2}} .. g_{\sigma_{n-2}} g_{\sigma_{n-1}} g_{\sigma_{n}} M'g_{\sigma_{n-1}} \bigg]  \nonumber\\\nonumber\\
				&=&\tau_{n+1} \bigg[ \underbrace{ g_{\sigma_{n-1}} g_{\sigma_{n}} g_{\sigma_{n-1}}}_{}  g_{\sigma_{n-2}} .. g_{\sigma_{2}} \big(^{n-1}_{2}F(^{2}g_{\sigma_{1}a_{2}})^{j} \big) g_{\sigma_{2}} ..  g_{\sigma_{n-2}} g_{\sigma_{n-1}} g_{\sigma_{n}} M' \bigg].\nonumber
			\end{eqnarray}
			
			We apply the TL relations. The cases corresponding to 1 and $g_{\sigma_{n-1}}$ are obvious.\\ 

			For the terms corresponding to $g_{\sigma_{n-1}}g_{\sigma_{n}} $, we have:
			\begin{eqnarray}
				& & \tau_{n+1} \bigg[ g_{\sigma_{n-1}} g_{\sigma_{n}} \underbrace{ g_{\sigma_{n-2}} .. g_{\sigma_{2}} \big(^{n-1}_{2} F(^{2}g_{\sigma_{1}a_{2}})^{j} \big) g_{\sigma_{2}} .. g_{\sigma_{n-2}}}_{} g_{\sigma_{n-1}} g_{\sigma_{n}} M' \bigg]  \nonumber\\\nonumber\\
				&=& \tau_{n+1} \bigg[ g_{\sigma_{n-1}} g_{\sigma_{n-2}} .. g_{\sigma_{2}} \big(^{n-1}_{2} F(^{2}g_{\sigma_{1}a_{2}})^{j} \big) g_{\sigma_{2}} .. g_{\sigma_{n-2}} \underbrace{ g_{\sigma_{n}} g_{\sigma_{n-1}} g_{\sigma_{n}}}_{} M' \bigg]. \nonumber
			\end{eqnarray}
			
			We are done, since it is a sum of values of $\tau_{n+1} $ on Markov elements, (the same for the term corresponding to $g_{\sigma_{n}} $) .  \\

			For the terms corresponding to  $g_{\sigma_{n}}g_{\sigma_{n-1}} $, we have:
			\begin{eqnarray}
				\tau_{n+1} \bigg[ g_{\sigma_{n}} g_{\sigma_{n-1}} g_{\sigma_{n-2}} .. g_{\sigma_{2}} \big(^{n-1}_{2}F(^{2} g_{\sigma_{1}a_{2}})^{j} \big) g_{\sigma_{2}} .. g_{\sigma_{n-2}} g_{\sigma_{n-1}} g_{\sigma_{n}} M' \bigg], \nonumber
			\end{eqnarray}
			
			which is the case of term (A), since $ M'$ is a Markov element in $\widehat{TL}_{n-1}(q) $.\\ 

			For terms of \textbf{Type (C)}, we see that:
			\begin{eqnarray}
				I_{3} &:=& \tau_{n+1} \bigg[ g_{\sigma_{n}} \big(g_{\sigma_{n-1} \sigma_{n-2} .. \sigma_{1}}  F_{n}(t_{a_{n}}) \big)^{k} g_{\sigma_{n}} F_{n} \big( F_{n-1} (e_{x_{i}}) (t_{\sigma_{n-1} \sigma_{n-2} ..\sigma_{1} a_{n}})^{i} \big) \bigg] \nonumber\\\nonumber\\
				&=& \tau_{n+1} \bigg[ g_{\sigma_{n}} F_{n} \big((t_{\sigma_{n-1} \sigma_{n-2} .. \sigma_{1}a_{n}} )^{k} \big) g_{\sigma_{n}} F_{n} \big(F_{n-1}(e_{x_{i}}) (t_{\sigma_{n-1} \sigma_{n-2} ..\sigma_{1} a_{n}})^{i} \big) \bigg]. \nonumber			
			\end{eqnarray}
			
			Call $ \gamma $ the image of $ F_{n-1}(e_{x_{i}}) $ under the action of $(t_{\sigma_{n-1} \sigma_{n-2} ..\sigma_{1} a_{n}})^{i}$. Thus: 
			\begin{eqnarray}
				I_{3} = \tau_{n+1} \bigg[ g_{\sigma_{n}} F_{n} \big((t_{\sigma_{n-1} \sigma_{n-2} .. \sigma_{1}a_{n}} )^{k} \big) g_{\sigma_{n}} F_{n} \big((t_{\sigma_{n-1} \sigma_{n-2} ..\sigma_{1} a_{n}})^{i} \gamma \big) \bigg]. \nonumber
			\end{eqnarray}
			
			As we have seen in the case (B), $(t_{\sigma_{n-1} \sigma_{n-2} .. \sigma_{1}a_{n}} )^{k} $ can be written as a sum of elements of the form (up to Markov elements in $\widehat{TL}_{n}(q) $ ):
			\begin{eqnarray}
				t_{\sigma_{n-1}} t_{\sigma_{n-2}} .. t_{\sigma_{2}} \big(^{n-1}_{2}F(^{2}g_{\sigma_{1}a_{2}})^{j} \big) t_{\sigma_{2}} .. t_{\sigma_{n-2}} t_{\sigma_{n-1}} \lambda, \nonumber
			\end{eqnarray} 
			
			where $ j \leq k $, and  $\lambda $ is in $^{n}_{n-1}F\big(\widehat{TL}_{n-1}(q)\big)$. \\

			Call $ \eta $ the image of $\lambda $ under the action of $(t_{\sigma_{n-1} \sigma_{n-2} ..\sigma_{1} a_{n}})^{i}$.\\

			The determination of $I_{3}$ can be reduced to computing the following value:
			\begin{eqnarray}
				\tau_{n+1} \bigg[ g_{\sigma_{n}}g_{\sigma_{n-1}} g_{\sigma_{n-2}} .. g_{\sigma_{2}} \big(^{n-1}_{2}F(^{2}g_{\sigma_{1}a_{2}})^{j} \big) g_{\sigma_{2}} .. g_{\sigma_{n-2}} g_{\sigma_{n-1}} g_{\sigma_{n}} F_{n}\big(\big(t_{\sigma_{n-1} \sigma_{n-2} ..\sigma_{1} a_{n}} \big)^{i} \eta \gamma \big)\bigg]. \nonumber
			\end{eqnarray}
			
			We repeat the same algorithm to $ (t_{\sigma_{n-1} \sigma_{n-2} ..\sigma_{1} a_{n}})^{i} $. Hence, we get some $l \leq i $, and some $\Delta$ in $^{n}_{n-1}F\big(\widehat{TL}_{n-1}(q)\big)$, such that we are reduced to compute:
			\begin{eqnarray}
				\tau_{n+1} &\bigg[& g_{\sigma_{n}}g_{\sigma_{n-1}}  g_{\sigma_{n-2}} ..  g_{\sigma_{2}} \big(^{n-1}_{2}F(^{2}g_{\sigma_{1}a_{2}})^{j} \big)  g_{\sigma_{2}} ..  g_{\sigma_{n-2}}  \underbrace{ g_{\sigma_{n-1}} g_{\sigma_{n}}g_{\sigma_{n-1}}}_{}  g_{\sigma_{n-2}} ..  g_{\sigma_{2}} \nonumber\\\nonumber\\
				& & \big(^{n-1}_{2}F(^{2}g_{\sigma_{1}a_{2}})^{l}\big)  g_{\sigma_{2}} ..  g_{\sigma_{n-2}}  g_{\sigma_{n-1}} \Delta \bigg]. \nonumber
			\end{eqnarray}

			We see, after using the T-L relations, that the terms corresponding to 1 and $ g_{\sigma_{n-1}}$ are values of $\tau_{n+1} $ on Markov elements.\\ 

			The term corresponding to $ g_{\sigma_{n-1}} g_{\sigma_{n}}$ is:
			\begin{eqnarray}
				\tau_{n+1} &\bigg[& g_{\sigma_{n}}g_{\sigma_{n-1}}  g_{\sigma_{n-2}} ..  g_{\sigma_{2}} \big(^{n-1}_{2}F(^{2}g_{\sigma_{1}a_{2}})^{j}\big)  g_{\sigma_{2}} ..  g_{\sigma_{n-2}}  g_{\sigma_{n-1}}  g_{\sigma_{n}}  g_{\sigma_{n-2}} ..  g_{\sigma_{2}} \nonumber\\\nonumber\\
				& & \big(^{n-1}_{2}F(^{2}g_{\sigma_{1}a_{2}})^{l}\big)  g_{\sigma_{2}} ..  g_{\sigma_{n-2}}  g_{\sigma_{n-1}} \Delta \bigg]  \nonumber\\\nonumber\\
				=\tau_{n+1} &\bigg[& g_{\sigma_{n-1}}  g_{\sigma_{n-2}} ..  g_{\sigma_{2}} \big(^{n-1}_{2}F(^{2}g_{\sigma_{1}a_{2}})^{j}\big)  g_{\sigma_{2}} ..  g_{\sigma_{n-2}}  g_{\sigma_{n-1}}   g_{\sigma_{n-2}} ..  g_{\sigma_{2}} \nonumber\\\nonumber\\
				&& \big(^{n-1}_{2}F(^{2}g_{\sigma_{1}a_{2}})^{l}\big)  g_{\sigma_{2}} ..  g_{\sigma_{n-2}} \underbrace{ g_{\sigma_{n}} g_{\sigma_{n-1}}g_{\sigma_{n}}}_{} \Delta \bigg]. \nonumber			
			\end{eqnarray}

			The term in square brackets is clearly a Markov element (the same thing with the term corresponding to $g_{\sigma_{n}}$) .\\

			The term corresponding to $ g_{\sigma_{n}} g_{\sigma_{n-1}} $ is:
			\begin{eqnarray}
				\tau_{n+1} &\bigg[& g_{\sigma_{n}}g_{\sigma_{n-1}}  g_{\sigma_{n-2}} ..  g_{\sigma_{2}} \big(^{n-1}_{2}F(^{2}g_{\sigma_{1}a_{2}})^{j} \big)  g_{\sigma_{2}} ..  g_{\sigma_{n-2}}  g_{\sigma_{n}}  g_{\sigma_{n-1}}  g_{\sigma_{n-2}} ..  g_{\sigma_{2}} \nonumber\\\nonumber\\
				&\big(&^{n-1}_{2}F(^{2}g_{\sigma_{1}a_{2}})^{l} \big)  g_{\sigma_{2}} ..  g_{\sigma_{n-2}}  g_{\sigma_{n-1}} \Delta \bigg]  \nonumber\\\nonumber\\
				=\tau_{n+1} &\bigg[& \underbrace{g_{\sigma_{n}}g_{\sigma_{n-1}} g_{\sigma_{n}}}_{} g_{\sigma_{n-2}} ..  g_{\sigma_{2}} \big(^{n-1}_{2}F(^{2}g_{\sigma_{1}a_{2}})^{j} \big)  g_{\sigma_{2}} ..  g_{\sigma_{n-2}}    g_{\sigma_{n-1}}  g_{\sigma_{n-2}} ..  g_{\sigma_{2}} \nonumber\\\nonumber\\
				&\big(&^{n-1}_{2}F(^{2}g_{\sigma_{1}a_{2}})^{l} \big)  g_{\sigma_{2}} ..  g_{\sigma_{n-2}}  g_{\sigma_{n-1}} \Delta \bigg]. \nonumber
			\end{eqnarray} 

			It is a Markov element, Theorem  \ref{5_1_1} follows.\\

	\section{Affine Markov trace: existence and uniqueness}

		\subsection{Existence}	
			Now, consider the following commutative diagram, where $(\tau_{n+1})_{0 \leq n}$ is the Markov trace from Theorem  \ref{1_1} and the vertical arrows $E_{n}: \widehat{TL}_{n+1}(q)  \longrightarrow  TL_{n}(q)$  have been defined in Proposition \ref{En}: \\
			
			\begin{tikzpicture}

			\matrix[matrix of math nodes,row sep=1.3cm,column sep=0.5cm]{
			|(A)|\widehat{TL}_{1}(q) & & & & |(B)| \widehat{TL}_{2}(q) & & |(BB)| \textcolor{white}{T}\dots ~~~~ & &|(C)|\widehat{TL}_{n}(q) & & & & |(D)| \widehat{TL}_{n+1}(q) \\
			|(E)|          TL_{0}(q) & & & & |(F)|           TL_{1}(q) & & |(FF)| \textcolor{white}{T}\dots ~~~~ & &|(G)| TL_{n-1}(q)        & & & & |(H)|             TL_{n}(q) \\
				                    & & & &                           & &                                       & &                         & & & &                             \\
				                    & & & &                           & &|(I) |K                                & &                         & & & &                             \\
				};

			\path (A) edge[-myto,line width=0.42pt]                                                                             (B);

			\path (B) edge[-myto,line width=0.42pt]                                                                             (BB);

			\path (BB) edge[-myto,line width=0.42pt]                                                                             (C);
			
			\path (C) edge[-myto,line width=0.42pt]                                                                             (D);
			
			
			\path (F) edge[-myhook,line width=0.42pt]                                                                           (E);
			\path (E) edge[-myto,line width=0.42pt]                                                                             (F);

			\path (FF) edge[-myhook,line width=0.42pt]                                                                           (F);
			\path (F) edge[-myto,line width=0.42pt]                                                                             (FF);
			
			\path (G) edge[-myhook,line width=0.42pt]                                                                           (FF);
			\path (FF) edge[-myto,line width=0.42pt]                                                                             (G);

			\path (H) edge[-myhook,line width=0.42pt]                                                                           (G);
			\path (G) edge[-myto,line width=0.42pt]                                                                             (H);


			\path (A) edge[-myonto,line width=0.42pt]                                                                           (E);
			
			\path (B) edge[-myonto,line width=0.42pt]                                                                           (F);

			\path (C) edge[-myonto,line width=0.42pt]                                                                           (G);

			\path (D) edge[-myonto,line width=0.42pt]                                                                           (H);

			\path (E) edge[-myto,line width=0.42pt]  node[above, xshift=-5mm, yshift=-2mm, rotate=0] {\footnotesize $\tau_{1}$}  (I);

			\path (F) edge[-myto,line width=0.42pt]  node[above, xshift=-3mm, yshift=-2mm, rotate=0] {\footnotesize $\tau_{2}$}  (I);

			\path (G) edge[-myto,line width=0.42pt]  node[above, xshift=3mm, yshift=-2mm, rotate=0] {\footnotesize $\tau_{n}$}  (I);

			\path (H) edge[-myto,line width=0.42pt]  node[above, xshift=8mm, yshift=-2mm, rotate=0] {\footnotesize $\tau_{n+1}$}  (I);

		\end{tikzpicture}

			\begin{proposition} \label{5_2_4} 
			Set $\rho_{n+1}$ to be the trace over  $\widehat{TL}_{n+1} (q) $ induced by $\tau_{n+1}$ over  $TL_{n} (q) $ for $0 \leq n$.
				We have: \\
				
				\begin{itemize}[label=$\bullet$, font=\normalsize, font=\color{black}, leftmargin=2cm, parsep=0cm, itemsep=0.25cm, topsep=0cm]
					\item $\rho_{n+1}(F_{n}(h)T^{\pm1}_{\sigma_{n}}) =  \rho_{n}(h)$, for all $h \in  \widehat{TL}_{n} (q)$, where $1 \leq n$.
					\item $\rho_{i}$ is invariant under the action of $\phi_{i}$ for all $i$. 
				\end{itemize}  
				

			
			\end{proposition}			

			\begin{demo}
			
			We have: $\rho_{n+1} \big(F_{n}(h)T^{\pm1}_{\sigma_{n}}\big) $ equals $ \tau_{n+1}\bigg(E_{n}\big(F_{n}(h)\big) E_{n}\big(T^{\pm1}_{\sigma_{n}}\big)\bigg)$.

				\begin{eqnarray}
					\text{Hence, }\rho_{n+1}\big(F_{n}(h)T^{\pm1}_{\sigma_{n}}\big) = \tau_{n+1}\bigg(x_{n}\big(E_{n-1}(h)\big)T^{\pm1}_{\sigma_{n}}\bigg) =  \tau_{n}\big(E_{n-1}(h)\big) = \rho_{n}\big(h\big). \nonumber 
				\end{eqnarray}
		
				We made use of the fact that the following diagram commutes, together with the fact that $(\tau_{n})_{1 \leq n}$ is a Markov trace: \\
				
				\begin{tikzpicture}

			\matrix[matrix of math nodes,row sep=0.9cm,column sep=0.75cm]{
			   |(A)| \widehat{TL_{n}(q)} & & & &                    &                 &                 & & & & |(B)| \widehat{TL_{n+1}(q)}   \\
			                             & & & &                    &                 &                 & & & &                               \\	
			                             & & & &  |(C)| TL_{n-1}(q) &                 & |(D)| TL_{n}(q) & & & &                               \\
				                        & & & &                    &                 &                 & & & &                                \\
				                        & & & &                    &                 &                 & & & &                                \\
				                        & & & &                    &                 &                 & & & &                                \\
			                             & & & &                    &|(E)| K          &                 & & & &                                \\	
				};
				
			\path (A) edge[-myto,line width=0.42pt] node[above, xshift=2mm, yshift=0mm, rotate=0] {\footnotesize $F_{n}$}     (B);
			
			\path (A) edge[-myonto,line width=0.42pt] node[above, xshift=2mm, yshift=0mm, rotate=0] {\footnotesize $E_{n-1}$}     (C);	
				
			\path (A) edge[-myto,line width=0.42pt] node[above, xshift=-7mm, yshift=0mm, rotate=0] {\footnotesize $\rho_{n}$}     (E);	
				
			\path (B) edge[-myonto,line width=0.42pt] node[above, xshift=-2mm, yshift=0mm, rotate=0] {\footnotesize $E_{n}$}     (D);	
		
			\path (B) edge[-myto,line width=0.42pt] node[above, xshift=7mm, yshift=0mm, rotate=0] {\footnotesize $\rho_{n+1}$}     (E);
			
			\path (C) edge[-myto,line width=0.42pt] node[above, xshift=2mm, yshift=0mm, rotate=0] {\footnotesize $\tau_{n}$}     (E);	
	
			\path (D) edge[-myto,line width=0.42pt] node[above, xshift=-4mm, yshift=0mm, rotate=0] {\footnotesize $\tau_{n+1}$}     (E);

			\path (D) edge[-myhook,line width=0.42pt] node[above, xshift=-2mm, yshift=0mm, rotate=0] {\footnotesize $x_{n} $}     (C);
			\path (C) edge[-myto,line width=0.42pt] node[above, xshift=-2mm, yshift=0mm, rotate=0] {\footnotesize $ $}     (D);

		\end{tikzpicture}

				For the second statement, we show that $\rho_{n}\big(h\big)=\rho_{n}\big([h]\big)$, where $[h]$ is the image of $h$ under $\phi^{-1}_{n}$. So we start from $\rho_{n}\big(h\big) =\tau_{n} \big(E_{n-1}(h)\big) $. But since $\tau_{n}$ is the $n$-th Markov trace, we have $\tau_{n} \big(E_{n-1}(h)\big) = -\frac{\sqrt{q}}{1+q}\tau_{n+1} \big(x_{n}(E_{n-1}(h))\big) $, which is equal to $ -\frac{\sqrt{q}}{1+q}\tau_{n+1} \bigg(E_{n}\big(F_{n}(h)\big)\bigg)$, since the diagram  $T$  commutes, this term is equal to $ -\frac{\sqrt{q}}{1+q}\rho_{n+1}\big(F_{n}(h)\big)$, hence to:
				\begin{eqnarray}
					-\frac{\sqrt{q}}{1+q}\rho_{n+1}\big(g_{\sigma_{n} .. \sigma_{1} a_{n+1}} F_{n}(h)g^{-1}_{\sigma_{n} .. \sigma_{1} a_{n+1}}\big) = -\frac{\sqrt{q}}{1+q}\rho_{n+1}\big( F_{n}([h])\big). \nonumber 
				\end{eqnarray}
				
				Now we consider the same steps in the opposite direction, that is:
				\begin{eqnarray}
					-\frac{\sqrt{q}}{1+q}\rho_{n+1}\big( F_{n}([h])\big) =  -\frac{\sqrt{q}}{1+q}\tau_{n+1} \bigg(E_{n}\big(F_{n}([h])\big)\bigg) = \rho_{n}\big([h]\big).  \nonumber
				\end{eqnarray}
				
			\end{demo}
		
			\begin{corollaire} \label{5_2_5}
				With the above notations, in the sense of Definition \ref{5_2_1}: $(\rho_{i})_{1 \leq i}$ is an affine Markov trace over $\big( \widehat{TL}_{i} (q)\big)_{1 \leq i}$. 
			\end{corollaire}

	 \subsection{Uniqueness}

			Consider the following homomorphism:
			\begin{eqnarray}
				F_{2}: \widehat{TL}_{2}(q) &\longrightarrow& \widehat{TL}_{3}(q) \nonumber \\
				\textsl{g}_{\sigma_{1}} &\longmapsto& g_{\sigma_{1}} \nonumber \\
				\textsl{g}_{a_{2}} &\longmapsto& g_{\sigma_{2}} g_{a_{3}} g^{-1}_{\sigma_{2}}  \nonumber\\\nonumber
			\end{eqnarray}
where, for possible lack of injectivity (see the comments preceding Proposition 
  \ref{2_1}), we use slanted letters  $\textsl{g}$, $\textsl{f}$  in $\widehat{TL}_{2}(q)$ 
  while we use the usual style $g$, $f$ in $\widehat{TL}_{3}(q)$. \\

			We set $F:=F_{2}$ in order to simplify in what follows. $F$ can be expressed by the following form considering the "f" generators, we see that $F\big(\textsl{f}_{a_{2}}\big) = F\big(\frac{\textsl{g}_{a_{2}}+1}{q+1} \big)$, which is equal to $\frac{1}{q+1} g_{\sigma_{2}} g_{a_{3}} g^{-1}_{\sigma_{2}}+\frac{1}{q+1}$, hence to: 
			\begin{eqnarray}
				\frac{1}{q+1} \bigg[\big((q+1) f_{\sigma_{2}} -1\big) \big((q+1) f_{a_{3}} -1\big)\big(\frac{1}{q}\big((q+1) f_{\sigma_{2}} -1\big)+\frac{1-q}{q}\big)\bigg]+\frac{1}{q+1}. \nonumber
			\end{eqnarray}
			
			Thus, we see that:
			\begin{eqnarray}
				F : \widehat{TL}_{2}(q) &\longrightarrow& \widehat{TL}_{3}(q) \nonumber\\
				\textsl{f}_{1} &\longmapsto& f_{\sigma_{1}}\nonumber\\
				\textsl{f}_{a_{2}} &\longmapsto& -\frac{q+1}{q}f_{a_{3}\sigma_{2}}-(q+1)f_{\sigma_{2}a_{3}}+f_{\sigma_{2}}+f_{a_{3}}.\nonumber
			\end{eqnarray}

			Notice that $F(\textsl{f}_{a_{2}})f_{\sigma_{2}}F(\textsl{f}_{a_{2}}) = \delta F(\textsl{f}_{a_{2}})$, and $f_{\sigma_{2}}F(\textsl{f}_{a_{2}}) f_{\sigma_{2}}=\delta f_{\sigma_{2}}$. Since we are interested with viewing $F\big(\widehat{TL}_{2}(q)\big)$ in $\widehat{TL}_{3}(q)$, we will investigate in what follows, the elements $\big(F(\textsl{f}_{\sigma_{1}}\textsl{f}_{a_{2}})\big)^{k}$ and $\big(F(\textsl{f}_{a_{2}}\textsl{f}_{\sigma_{1}})\big)^{k}$, for $k$ a positive integer.

			\begin{eqnarray}
				\text{Set }x_{1}:= F(\textsl{f}_{\sigma_{1}}\textsl{f}_{a_{2}})=f_{\sigma_{1}}F(\textsl{f}_{a_{2}})=  -\frac{q+1}{q}f_{\sigma_{1}a_{3}\sigma_{2}}-(q+1)f_{\sigma_{1}\sigma_{2}a_{3}}+f_{\sigma_{1}\sigma_{2}}+f_{\sigma_{1}a_{3}}. \nonumber
			\end{eqnarray}
		
			And for $1 \leq i$, we set: 
			\begin{eqnarray}
				x_{i}&:=&\big(-1\big)^{i}~~~\big(\frac{q+1}{q}\big)^{i}~~~~~f^{i}_{\sigma_{1}a_{3}\sigma_{2}}~~~~~+\big(-1\big)^{i}\big(q+1\big)^{i}f^{i}_{\sigma_{1}\sigma_{2}a_{3}}\nonumber\\\nonumber\\
				&+& \big(-1\big)^{i-1}\big(\frac{q+1}{q}\big)^{i-1}f^{i-1}_{\sigma_{1}a_{3}\sigma_{2}}f_{\sigma_{1}a_{3}}~~+\big(-1\big)^{i-1}\big(q+1\big)^{i-1}f^{i-1}_{\sigma_{1}\sigma_{2}a_{3}}f_{\sigma_{1}\sigma_{2}}.\nonumber
			\end{eqnarray}
		
			Notice that $ x^{2}_{1} = 3\delta x_{1}+ x_{2}$. It is easy to show  that:
			\begin{eqnarray}
				x_{1}x_{i} = \delta^{2} x_{i-1}+ 2\delta x_{i} + x_{i+1},  \text{  for } 2 \leq i, \nonumber
			\end{eqnarray}
		
			thus, for $ 1 \leq k $, we have $ x^{k}_{1}=\sum\limits^{i=k-1}_{i=1} \gamma_{i}x_{i}+ x_{k}$, here $\gamma_{i}$ is a polynomial in $\delta$, for all $i$. \\ 
		
			Notice that $x_{1}x_{j} = x_{j} x_{1}$ for $ j=1,2$. For $j=1$ it is clear, while for $j=2$ we have $ x_{2}=  x^{2}_{1} -3\delta x_{1}$. Now suppose that $3 \leq j$. We have $x_{j} = x_{1}x_{j-1} -\delta^{2} x_{j-2} -2\delta x_{j-1} $, hence we see by induction on $j$, that $x_{1}x_{j} = x_{j} x_{1}$, for all $j$ . \\
		
			We define the $\mathbb{Q}$-linear map $\chi: \widehat{TL}_{3}(q) \longrightarrow \widehat{TL}_{3}(q)$ which sends 1 to 1, and for any $u=s_{1}s_{2} .. s_{r}$ reduced expression of any element $u$ in $W^{c}(\tilde{A}_{2})$, it sends $f_{u}$ to $f_{s_{r}s_{r-1} .. s_{1}}$, with $q$ sent to $\frac{1}{q}$. \\
		
			Set $z_{1}:=F( \textsl{f}_{a_{2}} \textsl{f}_{\sigma_{1}})$. Then
			\begin{eqnarray}
				z_{1} =F(\textsl{f}_{a_{2}})f_{\sigma_{1}}=  -\frac{q+1}{q}f_{a_{3}\sigma_{2}\sigma_{1}}-(q+1)f_{\sigma_{2}a_{3}\sigma_{1}}+f_{\sigma_{2}\sigma_{1}}+f_{a_{3}\sigma_{1}}. \nonumber
			\end{eqnarray}
			
			And for $1 \leq i$, we set
			\begin{eqnarray}
				z_{i}:&=&\big(-1\big)^{i}~~~\big(\frac{q+1}{q}\big)^{i}~~~f^{i}_{a_{3}\sigma_{2}\sigma_{1}}~~~~~~~~+\big(-1\big)^{i}\big(q+1\big)^{i}f^{i}_{\sigma_{2}a_{3}\sigma_{1}} \nonumber\\\nonumber\\
				&+& \big(-1\big)^{i-1}\big(\frac{q+1}{q}\big)^{i-1}f_{\sigma_{2}\sigma_{1}}f^{i-1}_{a_{3}\sigma_{2}\sigma_{1}}~~+\big(-1\big)^{i-1}\big(q+1\big)^{i-1}f_{a_{3}\sigma_{1}}f^{i-1}_{\sigma_{2}a_{3}\sigma_{1}}. \nonumber
			\end{eqnarray}
			
			Notice that $\chi(x_{i}) = z_{i}$ for all $i $. Now $\chi(x_{1}x_{i})=\chi(x_{i}x_{1})=z_{1}z_{i}$. We see that $\chi(\delta) = \delta$. Moreover, $z_{1}z_{j} = \chi (x_{1}x_{j}) = \chi (\delta^{2} x_{i-1}+ 2\delta x_{i} + x_{i+1}) = \delta^{2} z_{i-1}+ 2\delta z_{i} + z_{i+1}$. And in the same way, by acting by $\chi$, we find that $z^{k}_{1}=\sum\limits^{i=k-1}_{i=1} \gamma_{i}z_{i}+ z_{k}$, where $\gamma_{i}$ is as above.\\
		
		Consider $ x_{i}f_{\sigma_{2}}$ for $1 \leq i$, we see that it is equal to:\\
		\begin{eqnarray}
			& & \big(-1\big)^{i}~~~\big(\frac{q+1}{q}\big)^{i}~~~f^{i}_{\sigma_{1}a_{3}\sigma_{2}}~~~~~~~~~~~~+\big(-1\big)^{i}\big(q+1\big)^{i}f^{i}_{\sigma_{1}\sigma_{2}a_{3}}f_{\sigma_{2}}  \nonumber\\\nonumber\\
			&+& \big(-1\big)^{i-1}\big(\frac{q+1}{q}\big)^{i-1}f^{i-1}_{\sigma_{1}a_{3}\sigma_{2}}f_{\sigma_{1}a_{3}}f_{\sigma_{2}}~~+\big(-1\big)^{i-1}\big(q+1\big)^{i-1}f^{i-1}_{\sigma_{1}\sigma_{2}a_{3}}f_{\sigma_{1}\sigma_{2}}, \nonumber
		\end{eqnarray}

		\begin{eqnarray}
			\text{which is: }& & \big(-1\big)^{i}~\big(\frac{q+1}{q}\big)^{i}~~~f^{i}_{\sigma_{1}a_{3}\sigma_{2}}~~~~~~~~~~~+\delta \big(-1\big)^{i}\big(q+1\big)^{i}f^{i-1}_{\sigma_{1}\sigma_{2}a_{3}}f_{\sigma_{1}\sigma_{2}} \nonumber\\\nonumber\\
			&+& \big(-1\big)^{i}\big(\frac{q+1}{q}\big)^{i-1}f^{i}_{\sigma_{1}a_{3}\sigma_{2}}~~~~~~~~~~~~+\big(-1\big)^{i-1}\big(q+1\big)^{i-1}f^{i-1}_{\sigma_{1}\sigma_{2}a_{3}}f_{\sigma_{1}\sigma_{2}}.\nonumber
		\end{eqnarray}

		\begin{eqnarray}
			\text{Hence, }x_{i}f_{\sigma_{2}} = \big[(-1)^{i-1}\frac{(q+1)^{i-1}}{q^{i}} \big] f^{i}_{\sigma_{1}a_{3}\sigma_{2}} +  \big[(-1)^{i-1}(q+1)^{i-2} \big] f^{i-1}_{\sigma_{1}\sigma_{2}a_{3}}f_{\sigma_{1}\sigma_{2}},\text{ for } 1 \leq i. \nonumber
		\end{eqnarray}
		
		\vspace{-0.5 cm}
		
		\begin{eqnarray}
		    \text{In particular }x_{1}f_{\sigma_{2}} = -\frac{q+1}{q}f_{\sigma_{1}a_{3}\sigma_{2}}-(q+1)f_{\sigma_{1}\sigma_{2}a_{3}}f_{\sigma_{2}}+f_{\sigma_{1}\sigma_{2}}+f_{\sigma_{1}a_{3}}f_{\sigma_{2}}, \nonumber
		\end{eqnarray}
				
		thus,
		\vspace{-0.5 cm}
		
		\begin{eqnarray}
		 x_{1}f_{\sigma_{2}}= \frac{-1}{q}f_{\sigma_{1}a_{3}\sigma_{2}}+ \frac{1}{q+1}f_{\sigma_{1}\sigma_{2}}. \nonumber
		\end{eqnarray}

		Now we apply $\chi$ to $ x_{i}f_{\sigma_{2}}$. Hence
		\begin{eqnarray}
			f_{\sigma_{2}} z_{i} = \big[(-1)^{i}q(q+1)^{i-1}\big] f^{i}_{\sigma_{2}a_{3}\sigma_{1}} + \big[(-1)^{i-1}(\frac{q+1}{q})^{i-2}\big]f_{\sigma_{2}\sigma_{1}}f^{i-1}_{a_{3}\sigma_{2}\sigma_{1}},\text{ for } 1 \leq i. \nonumber
		\end{eqnarray}
		
		In particular $ f_{\sigma_{2}}z_{1}= -qf_{\sigma_{2}a_{3}\sigma_{1}}+\frac{q}{q+1}f_{\sigma_{2}\sigma_{1}}$.\\ 
		
		We consider Propositions \ref{5_1_10} and \ref{3_4}, take $t$ to be any $\psi_{2}$-invariant trace over $\widehat{TL}_{2} (q)$, determined by $A_{0},A_{1}$ and $(\alpha_{i})_{1\leq i}$. Let $s$ be any $\psi_{3}$-invariant trace over $\widehat{TL}_{3} (q)$, determined by $B_{0},B_{1},B_{2}$ and $(\beta_{i})_{1\leq i}$. We show in what follows that: \\ 
		
	\textit{ there are a unique $t$ and a unique $s$, such that $t$ is the second component and $s$ is the third component of a Markov trace.} 
		
		In order to simplify, we set $\hat{\tau}_{2}:=t$ and $\hat{\tau}_{3}:=s$.\\
		
		At first, being a first component of a Markov trace, forces $\hat{\tau}_{2}$ to have the value 1 on $\textbf{T}_{\sigma_{1}}$ and $\textbf{T}_{a_{2}}$, but $\textsl{f}_{\sigma_{1}} = \frac{1+\textsl{g}_{\sigma_{1}}}{1+q} = \frac{1}{1+q}+\frac{\textbf{T}_{\sigma_{1}}}{\sqrt{q}(1+q)}$. Hence, $A_{1} =- \frac{\sqrt{q}}{1+q}$. Moreover, $\hat{\tau}_{2}(1) = -\frac{1+q}{\sqrt{q}}\hat{\tau}_{1}(1) $. Thus, $A_{0}= - \frac{1+q}{\sqrt{q}}$.\\
		
		Now, we have: 
		\begin{eqnarray}
			B_{0} = \hat{\tau}_{3}\big(1\big) = - \frac{1+q}{\sqrt{q}}\hat{\tau}_{2}\big(1\big) = \big(- \frac{1+q}{\sqrt{q}}\big)^{2}, \nonumber
		\end{eqnarray}

		\begin{eqnarray}
			\text{and }B_{1} = \hat{\tau}_{3}\big(f_{\sigma_{1}}\big) =- \frac{1+q}{\sqrt{q}}\hat{\tau}_{2}\big(\textsl{f}_{\sigma_{1}}\big) = \frac{1+q}{\sqrt{q}}\frac{\sqrt{q}}{1+q} = 1.\nonumber
		\end{eqnarray}
		
		\begin{remarque}
			$\hat{\tau}_{3}$ must verify $\hat{\tau}_{3}\big( F(h)T_{\sigma_{2}}\big) =\hat{\tau}_{2}\big( h\big) $, for every $h$ in $\widehat{TL}_{2}(q)$. 
			\begin{eqnarray}
				\text{But, }\hat{\tau}_{3}\big( F(h)T_{\sigma_{2}}\big)=\sqrt{q}\hat{\tau}_{3}\big( F(h)g_{\sigma_{2}}\big)=\sqrt{q}\hat{\tau}_{3}\bigg( F(h)\big[(q+1)f_{\sigma_{2}} -1\big]\bigg). \nonumber
			\end{eqnarray}

			\begin{eqnarray}
				\text{So, }\sqrt{q}\hat{\tau}_{3}\bigg( F(h)\big[(q+1)f_{\sigma_{2}} -1\big]\bigg) &=&  \sqrt{q}\big(q+1\big)\hat{\tau}_{3}\big( F(h)f_{\sigma_{2}}\big)- \sqrt{q}\hat{\tau}_{3}\big( F(h)\big) \nonumber\\\nonumber\\
				& & \sqrt{q}\big(q+1\big)\hat{\tau}_{3}\big( F(h)f_{\sigma_{2}}\big)+ \sqrt{q}\frac{1+q}{\sqrt{q}}\hat{\tau}_{2}\big( h\big) .\nonumber				
			\end{eqnarray}
			
			Hence, our condition becomes
			\begin{eqnarray}
				\sqrt{q}\big(q+1\big)\hat{\tau}_{3}\big( F(h)f_{\sigma_{2}}\big) =  -\sqrt{q}\frac{1+q}{\sqrt{q}}\hat{\tau}_{2}\big( h\big) + \hat{\tau}_{2}\big( h\big) = -q\hat{\tau}_{2}\big( h\big). \nonumber
			\end{eqnarray}
	
			Thus, we must have
			\begin{eqnarray}
				\hat{\tau}_{3}\big( F(h)f_{\sigma_{2}}\big) = -\frac{\sqrt{q}}{(q+1)} \hat{\tau}_{2}\big( h\big), \text{ as an "f" equivalent to }\hat{\tau}_{3}\big( F(h)T_{\sigma_{2}}\big) =\hat{\tau}_{2}\big( h\big). \nonumber
			\end{eqnarray}
			
		\end{remarque}
		
	Now, we have:
		\vspace{-0.8 cm}
		
		\begin{eqnarray}
			B_{2} =  \hat{\tau}_{3}(f_{\sigma_{1}\sigma_{2}}) = -\frac{\sqrt{q}}{1+q}\hat{\tau}_{2}(\textsl{f}_{\sigma_{1}}) = ( \frac{\sqrt{q}}{1+q})^{2}.\nonumber
		\end{eqnarray}
		
		So, under the assumption that our two traces are the second and the third components of a given Markov trace, we get the following:

		\begin{eqnarray}
			A_{1} &=& -\frac{\sqrt{q}}{1+q}$, $A_{0}=  -\frac{1+q}{\sqrt{q}}. \nonumber\\\nonumber\\
			B_{2} &=& \big( \frac{\sqrt{q}}{1+q}\big)^{2} $, $ B_{1}= 1 \text{ and } B_{0}= \big( \frac{1+q}{\sqrt{q}}\big)^{2}.\nonumber
		\end{eqnarray}
		
		In particular, we have for all $1 \leq i$:
		\begin{eqnarray}
			\hat{\tau}_{3}\big( x^{i}_{1}f_{\sigma_{2}}\big) = -\frac{\sqrt{q}}{(q+1)} \hat{\tau}_{2}\big(( \textsl{f}_{\sigma_{1}a_{2}})^{i}\big), ~~~~\text{and}~~~~\hat{\tau}_{3}\big( z^{i}_{1}f_{\sigma_{2}}\big) = -\frac{\sqrt{q}}{(q+1)} \hat{\tau}_{2}\big(( \textsl{f}_{a_{2}\sigma_{1}})^{i}\big). \nonumber
		\end{eqnarray}
		
		In other terms, for all $i$ we have:
		\begin{eqnarray}
			\hat{\tau}_{3}(x^{i}_{1}f_{\sigma_{2}}) = -\frac{\sqrt{q}}{(q+1)}\alpha_{i},~~~~\text{and}~~~~\hat{\tau}_{3}(f_{\sigma_{2}}z^{i}_{1}) = -\frac{\sqrt{q}}{(q+1)}\alpha_{i}. \nonumber
		\end{eqnarray}
		
		Since $\hat{\tau}_{3}$ is determined by $\beta_{i}$, we can view these equalities as system of equations in  $\beta_{i}$ and $\alpha_{i}$. In what follows, we show that this system has at most one solution: $(\alpha_{i},\beta_{i})_{1 \leq i}$. \\
				
		For $i=1$, we see that we have two equations:
		\begin{eqnarray}
			& & \hat{\tau}_{3}(\frac{-1}{q}f_{\sigma_{1}a_{3}\sigma_{2}}+ \frac{1}{q+1}f_{\sigma_{1}\sigma_{2}}) = -\frac{\sqrt{q}}{(q+1)}\alpha_{1},~~~~\text{and}~~~~\hat{\tau}_{3}(-qf_{\sigma_{2}a_{3}\sigma_{1}}+\frac{q}{q+1}f_{\sigma_{2}\sigma_{1}}) = -\frac{\sqrt{q}}{(q+1)}\alpha_{1},  \nonumber\\ \nonumber\\
			& & \text{that is} ~~~\frac{-1}{q}\beta_{1}+ \frac{1}{q+1}B_{2} = -\frac{\sqrt{q}}{(q+1)}\alpha_{1,}~~~~\text{and}~~~~-q\beta_{1}+\frac{q}{q+1}B_{2} = -\frac{\sqrt{q}}{(q+1)}\alpha_{1}, \nonumber\\ \nonumber\\
			& & \text{that is} ~~~\frac{-1}{q}\beta_{1}+ \frac{q}{(q+1)^{3}} =- \frac{\sqrt{q}}{(q+1)}\alpha_{1},~~~~\text{and}~~~~-q\beta_{1}+\frac{q^{2}}{(q+1)^{3}} = -\frac{\sqrt{q}}{(q+1)}\alpha_{1}. \nonumber
		\end{eqnarray}
		
		Clearly, those two linear equations are independent, hence, they determine a unique solution $(\alpha_{1},\beta_{1})$. Let us see the equations when $i=2$, we have:
		\begin{eqnarray}
			\hat{\tau}_{3}( x^{2}_{1}f_{\sigma_{2}}) = -\frac{\sqrt{q}}{(q+1)}\alpha_{2},~~~~\text{and}~~~~\hat{\tau}_{3}(f_{\sigma_{2}} z^{2}_{1}) = -\frac{\sqrt{q}}{(q+1)}\alpha_{2}. \nonumber
		\end{eqnarray}
				
		We see that:
		\begin{eqnarray}
			x^{2}_{1}f_{\sigma_{2}} &=& 3\delta x_{1}f_{\sigma_{2}}+ x_{2}f_{\sigma_{2}}=3\frac{-1}{q}\delta f_{\sigma_{1}a_{3}\sigma_{2}}+ 3\frac{1}{q+1}\delta f_{\sigma_{1}\sigma_{2}} -\frac{(q+1)}{q^{2}} f^{2}_{\sigma_{1}a_{3}\sigma_{2}} - f_{\sigma_{1}\sigma_{2}a_{3}}f_{\sigma_{1}\sigma_{2}}\nonumber\\ \nonumber\\
			&=& \frac{-3}{(1+q)^{2}} f_{\sigma_{1}a_{3}\sigma_{2}}+ \frac{3}{(1+q)^{3}} f_{\sigma_{1}\sigma_{2}} -\frac{(q+1)}{q^{2}} f^{2}_{\sigma_{1}a_{3}\sigma_{2}} - f_{\sigma_{1}\sigma_{2}a_{3}\sigma_{1}\sigma_{2}},\nonumber
		\end{eqnarray}

		\begin{eqnarray}
			\text{hence, } \hat{\tau}_{3}( x^{2}_{1}f_{\sigma_{2}}) &=& \frac{-3}{(1+q)^{2}}\beta_{1}+ \frac{3}{(1+q)^{3}} B_{2} -\frac{(q+1)}{q^{2}} \beta_{2} - \delta \beta_{1} \nonumber\\\nonumber\\
			&=& \frac{3}{(1+q)^{3}} B_{2}- \frac{3+q}{(1+q)^{2}}\beta_{1} -\frac{(q+1)}{q^{2}} \beta_{2}. \nonumber
		\end{eqnarray}
	
		\begin{eqnarray}
			\text{Now, }f_{\sigma_{2}} z^{2}_{1}&=& \chi( x^{2}_{1}f_{\sigma_{2}}) \nonumber\\\nonumber\\
			&=& \chi \bigg(\frac{-3}{(1+q)^{2}} \bigg)f_{\sigma_{2}a_{3}\sigma_{1}}+ \chi\bigg(\frac{3}{(1+q)^{3}} \bigg)f_{\sigma_{2}\sigma_{1}} -\chi\bigg(\frac{(q+1)}{q^{2}}\bigg) f^{2}_{\sigma_{2}a_{3}\sigma_{1}} - f_{\sigma_{2}\sigma_{1}a_{3}\sigma_{2}\sigma_{1}}, \nonumber
		\end{eqnarray}
		
		\begin{eqnarray}
			\text{so }f_{\sigma_{2}} z^{2}_{1} = \frac{-3q^{2}}{(1+q)^{2}}f_{\sigma_{2}a_{3}\sigma_{1}}+\frac{3q^{3}}{(1+q)^{3}}f_{\sigma_{2}\sigma_{1}} -q(q+1)f^{2}_{\sigma_{2}a_{3}\sigma_{1}} - f_{\sigma_{2}\sigma_{1}a_{3}\sigma_{2}\sigma_{1}}. \nonumber
		\end{eqnarray}

		\begin{eqnarray}
		\text{Now, we apply the trace }	\hat{\tau}_{3}\big(f_{\sigma_{2}} z^{2}_{1}\big) &=& \frac{-3q^{2}}{(1+q)^{2}}\beta_{1}+\frac{3q^{3}}{(1+q)^{3}}B_{2} -q(q+1)\beta_{2} - \delta \beta_{1}  \nonumber\\\nonumber\\
			&=& \frac{3q^{3}}{(1+q)^{3}}B_{2}-\frac{3q^{2}+q}{(1+q)^{2}}\beta_{1} -q(q+1)\beta_{2}. \nonumber
		\end{eqnarray}
		
		In other terms, we have the two equations:\\
		\begin{eqnarray}
			-\frac{(q+1)}{q^{2}} \beta_{2}- \frac{3+q}{(1+q)^{2}}\beta_{1} +\frac{3q}{(1+q)^{5}} &=& -\frac{\sqrt{q}}{(q+1)}\alpha_{2}, \nonumber\\\nonumber\\			
			-q(q+1)\beta_{2} - \frac{3q^{2}+q}{(1+q)^{2}}\beta_{1}+\frac{3q^{4}}{(1+q)^{5}} &=& - \frac{\sqrt{q}}{(q+1)}\alpha_{2},\nonumber
		\end{eqnarray}
	which indeed determine a unique  $(\alpha_{2},\beta_{2})$ as a solution.\\
		
		Now, we have: 
		\begin{eqnarray}
			x^{k}_{1}=\sum^{i=k-1}_{i=1} \gamma_{i}x_{i}+ x_{k},\nonumber
		\end{eqnarray}

		\begin{eqnarray}
			\text{hence, }x^{k}_{1}f_{\sigma_{2}}=\sum^{i=k-1}_{i=1} \gamma_{i}x_{i}f_{\sigma_{2}}+ x_{k}f_{\sigma_{2}},\nonumber
		\end{eqnarray}
		
		thus
		\begin{eqnarray}
			x^{k}_{1}f_{\sigma_{2}} = \sum^{i=k-1}_{i=1} \gamma_{i}&\bigg[&\big(-1\big)^{i-1}\frac{(q+1)^{i-1}}{q^{i}}\bigg] f^{i}_{\sigma_{1}a_{3}\sigma_{2}} +\gamma_{i} \bigg[\big(-1\big)^{i-1}(q+1)^{i-2}\bigg] f^{i-1}_{\sigma_{1}\sigma_{2}a_{3}}f_{\sigma_{1}\sigma_{2}} \nonumber\\\nonumber\\\nonumber\\
			 + &\bigg[&\big(-1\big)^{k-1}\frac{(q+1)^{k-1}}{q^{k}}\bigg] f^{k}_{\sigma_{1}a_{3}\sigma_{2}} + \bigg[\big(-1\big)^{k-1}(q+1)^{k-2}\bigg] f^{k-1}_{\sigma_{1}\sigma_{2}a_{3}}f_{\sigma_{1}\sigma_{2}}. \nonumber			
		\end{eqnarray}
		
		Now we apply $\hat{\tau}_{3}$, we get:
		\begin{eqnarray}
			-\frac{\sqrt{q}}{(q+1)} \alpha_{k}= \sum^{i=k-1}_{i=1} \gamma_{i}&\bigg[&\big(-1\big)^{i-1}\frac{(q+1)^{i-1}}{q^{i}}\bigg] \beta_{i} +\delta\gamma_{i} \bigg[\big(-1\big)^{i-1}(q+1)^{i-2}\bigg] \beta_{i-1}\nonumber\\\nonumber\\\nonumber\\
			+ \delta &\bigg[&\big(-1\big)^{k-1}(q+1)^{k-2}\bigg] \beta_{k-1}+\bigg[\big(-1\big)^{k-1}\frac{(q+1)^{k-1}}{q^{k}}\bigg] \beta_{k}. \nonumber			
		\end{eqnarray}
		
		It is clear that the coefficient of $\beta_{k}$ is not zero, since $\beta_{k}$ does not appear in:
		\begin{eqnarray}
			A:= \sum^{i=k-1}_{i=1} \gamma_{i}&\bigg[&\big(-1\big)^{i-1}\frac{(q+1)^{i-1}}{q^{i}}\bigg] \beta_{i} +\delta\gamma_{i} \bigg[\big(-1\big)^{i-1}(q+1)^{i-2}\bigg] \beta_{i-1} \nonumber\\\nonumber\\\nonumber\\
			+ \delta &\bigg[&\big(-1\big)^{k-1}(q+1)^{k-2}\bigg] \beta_{k-1}.\nonumber
		\end{eqnarray}
		
		Now, we repeat the same steps with $z_{i}$, namely:
		\begin{eqnarray}
			z^{k}_{1}=\sum^{i=k-1}_{i=1} \gamma_{i}d_{i}+ d_{k}, \nonumber
		\end{eqnarray}

		\begin{eqnarray}
			\text{hence, }f_{\sigma{2}}z^{k}_{1}=\sum^{i=k-1}_{i=1} \gamma_{i}f_{\sigma_{2}}d_{i}+ f_{\sigma_{2}}d_{k}.\nonumber
		\end{eqnarray}
		
		Thus,
		\begin{eqnarray}
			f_{\sigma_{2}}z^{k}_{1} = \sum^{i=k-1}_{i=1} \gamma_{i} &\bigg[&\big(-1\big)^{i}q(q+1)^{i-1}\bigg] f^{i}_{\sigma_{2}a_{3}\sigma_{1}} + \gamma_{i} \bigg[\big(-1\big)^{i-1}(\frac{q+1}{q})^{i-2}\bigg]f_{\sigma_{2}\sigma_{1}}f^{i-1}_{a_{3}\sigma_{2}\sigma_{1}} \nonumber\\\nonumber\\\nonumber\\
			+ &\bigg[&\big(-1\big)^{k}q(q+1)^{k-1}\bigg] f^{k}_{\sigma_{2}a_{3}\sigma_{1}} + \bigg[\big(-1\big)^{k-1}(\frac{q+1}{q})^{k-2}\bigg]f_{\sigma_{2}\sigma_{1}}f^{k-1}_{a_{3}\sigma_{2}\sigma_{1}}. \nonumber			
		\end{eqnarray}
		
		Now we apply $\hat{\tau}_{3}$, we get:
		\begin{eqnarray}
			-\frac{\sqrt{q}}{(q+1)} \alpha_{k} = \sum^{i=k-1}_{i=1} \gamma_{i} &\bigg[&\big(-1\big)^{i}q(q+1)^{i-1}\bigg] \beta_{i} + \gamma_{i} \delta\bigg[\big(-1\big)^{i-1}(\frac{q+1}{q})^{i-2}\big]\beta_{i-1} \nonumber\\\nonumber\\\nonumber\\		
		+ \delta &\bigg[&\big(-1)^{k-1}(\frac{q+1}{q})^{k-2}\bigg]\beta_{k-1}+\bigg[\big(-1\big)^{k}q(q+1)^{k-1}\bigg]\beta_{k}. \nonumber
		\end{eqnarray}
		
		The coefficient of $\beta_{k}$ is not zero, since $\beta_{k}$ does not appear in 
		\begin{eqnarray}
			B:=  \sum^{i=k-1}_{i=1} \gamma_{i} &\bigg[&\big(-1\big)^{i}q(q+1)^{i-1}\bigg] \beta_{i} + \gamma_{i} \delta\bigg[\big(-1\big)^{i-1}(\frac{q+1}{q})^{i-2}\bigg]\beta_{i-1} \nonumber\\\nonumber\\\nonumber\\	
			+ \delta&\bigg[&\big(-1)^{k-1}(\frac{q+1}{q})^{k-2}\bigg]\beta_{k-1}. \nonumber
		\end{eqnarray}
		
		In other terms, we have the two following equations, in $\beta_{k}$ and $\alpha_{k} $:\\
		\begin{eqnarray}
			-\frac{\sqrt{q}}{(q+1)} \alpha_{k} &=& A+\bigg[\big(-1\big)^{k-1}\frac{(q+1)^{k-1}}{q^{k}}\bigg] \beta_{k},\nonumber\\\nonumber\\\nonumber\\				
			-\frac{\sqrt{q}}{(q+1)} \alpha_{k}&=& B+\bigg[\big(-1\big)^{k}q(q+1)^{k-1}\bigg]\beta_{k}. \nonumber
		\end{eqnarray}
		
		Those are two independent linear equations in $\beta_{k}$ and $\alpha_{k}$, with non-zero coefficients, by induction over $k$ (that is: assuming that  $(\alpha_{i},\beta_{i})$ is unique for $ i < k$ then $(\alpha_{k},\beta_{k})$ is unique) we get the following corollary.\\ 
		
		\begin{corollaire} \label{5_2_7} 
			Suppose that $(\hat{\tau}_{i})_{1 \leq i}$ is a Markov trace over the tower of $\tilde{A}$-type T-L algebras, then $\hat{\tau}_{i} = \rho_{i}$ for $i = 1,2,3$.	\\	
		\end{corollaire}
		  
		Finally, we sum up the proof of the main theorem: we know, by Corollary \ref{5_2_5},  that there exists, at least, one affine Markov trace. Now, Corollary \ref{5_2_7} says that in any given affine Markov trace, the three first components are $\rho_{1},\rho_{2}$ and $\rho_{3}$ (of Corollary \ref{5_2_5}), while Theorem \ref{5_1_1} affirms that a third component in a given Markov trace determines a unique fourth component, and so on for any $\hat{\tau}_{i}$ with $i\geq 3$. Hence, we get our main theorem:  \\

		\begin{theoreme}\label{5_2_8} 
		
	 {	
		There exists a unique affine Markov trace over the tower of $\tilde{A}$-type Temperley-Lieb algebras, namely $(\rho_{i})_{1 \leq i}$. }\\
	\end{theoreme}

%
	

\renewcommand{\refname}{REFERENCES}

\end{document}